\documentclass[11pt]{article}
\usepackage{latexsym}
\usepackage{amsfonts,amssymb,amsmath}

\usepackage[yyyymmdd]{datetime}

\newtheorem{thm}{Theorem}[section]
\newtheorem{cor}[thm]{Corollary}
\newtheorem{lem}[thm]{Lemma}
\newtheorem{pro}[thm]{Proposition}
\newtheorem{defn}[thm]{Definition}

\title{Inverse monoids associated with the complexity class {\sf NP}}

\author{ J.C.\ Birget  }

\date{\footnotesize{6 March 2017}}
\begin{document}
\maketitle

\begin{abstract}
We study the {\sf P} versus {\sf NP} problem through properties of functions
and monoids, continuing the work of \cite{s1f}. Here we consider inverse 
monoids whose properties and relationships determine whether {\sf P} is 
different from {\sf NP}, or whether injective one-way functions (with respect 
to worst-case complexity) exist. 
\end{abstract}

%%%%%%%%%%%%%%%%%%%%%
\section{Introduction}

We give a few definitions before motivating the monoid approach to the 
{\sf P} versus {\sf NP} problem. Some of these notions appeared already in 
\cite{s1f}. 

\smallskip

By {\em function} we mean a {\em partial} function $A^* \to A^*$, where 
$A$ is a finite alphabet (usually, $A = \{0,1\}$), and $A^*$ denotes the 
set of all finite strings over $A$.
Let ${\sf Dom}(f)$ ($\subseteq A^*$) denote the domain of $f$, i.e.,
$\{x \in A^* : f(x)$ is defined\}; and let ${\sf Im}(f)$ ($\subseteq A^*$)
denote the image (or range) of  $f$, i.e., $\{f(x) : x \in {\sf Dom}(f)\}$.
The length of $x \in A^*$ is denoted by $|x|$.
The restriction of $f$ to $X \subseteq A^*$ is denoted by $f|_X$, and the 
identity function on $X$ is denoted by ${\sf id}_X$ or ${\sf id}|_X$.

\begin{defn} \label{InvCoMutinv} {\bf (inverse, co-inverse, mutual
inverse)}  
 \ A function $f'$ is an {\em inverse} of a function $f$ iff 
$\, f \circ f' \circ f = f$. In that case we also say that $f$ is a 
{\em co-inverse} of $f'$. \ If $f'$ is both an inverse and a co-inverse of 
$f$, we say that $f'$ is a {\em mutual inverse} of $f$, or that $f'$ and 
$f$ are mutual inverses (of each other).\footnote{The terminology about 
 {\it inverses} varies. 
 In semigroup theory, $f'$ such that $f f' f = f$ is called a semi-inverse 
 or a pseudo-inverse of $f$, in numerical mathematics $f'$ is called a  
 generalized inverse, in ring theory and in category theory it is called a 
 weak inverse or a von Neumann inverse. In semigroup theory the term 
 ``inverse'' of $f$ is only used if $f' f f' = f'$ holds in addition to 
 $f f' f = f$.  Co-inverses don't seem to have a name in the literature. 
}
\end{defn}
It is easy to see that $f'$ is an inverse of $f$ iff for every 
$y \in {\sf Im}(f)$: $f'(y)$ is defined and $f'(y) \in f^{-1}(y)$.
Hence, in particular, ${\sf Im}(f) \subseteq {\sf Dom}(f')$.
Also,  $f'$ is an inverse of $f$ iff
$\, f \circ f'|_{{\sf Im}(f)} = {\sf id}_{{\sf Im}(f)}$.

An element $f$ of a semigroup $S$ is called {\em regular} (or more
precisely, regular in $S$) iff there exists $f' \in S$ such that 
$f f' f = f$; in that case, $f'$ is called an {\em inverse} of $f$ (more
precisely, an inverse in $S$ of $f$), and $f$ is called a {\em co-inverse} 
of $f'$ (in $S$).
A semigroup $S$ is called regular iff all the elements of $S$ are regular 
in $S$; see e.g.\ \cite{Grillet,CliffPres}.
An {\em inverse semigroup} is, by definition, a semigroup $S$ such that every
element of $S$ has exactly {\em one} mutual inverse in $S$; 
see e.g.\ \cite{LawsonInv, Grillet,CliffPres}.

\medskip

The connection between the {\sf P} versus {\sf NP} problem and inverses 
comes from the following well known characterization:  
$\, {\sf P} \neq {\sf NP} \,$ iff there exists a function that is computable 
in polynomial time (by a deterministic Turing machine) and that is
polynomially balanced, but that does not have an inverse that is computable 
in polynomial time; see e.g.\ \cite{HemaOgi} p.\ 33, and Section 1 of 
\cite{s1f}. 
A function $f$: $A^* \to A^*$ is called {\em polynomially balanced} iff
there exists a polynomial $p$ such that for all $x \in {\sf Dom}(f)$:
 \, $|f(x)| \leq p(|x|)$ and $|x| \leq p(|f(x)|)$.
Functions that are polynomial-time computable and polynomially balanced, but
that have no inverse of that type, are {\em one-way} functions in the sense of worst-case 
 complexity.\footnote{One-way functions in this sense are not necessarily
 useful for cryptography, but they are relevant for {\sf P} vs.\ {\sf NP}.
 }  %% [End, footnote]
In \cite{s1f} this characterization of the {\sf P} versus {\sf NP} problem 
was reformulated in terms of regularity of semigroups, by introducing the 
monoid {\sf fP}, defined as follows.

\begin{defn} \label{DEFfP} \hspace{-0.07in}{\bf .} 

${\sf fP} \ = \ \{f : f$ is a function $A^* \to A^*$ that is polynomially
balanced, and 

\hspace{.7in} 
computable by a deterministic polynomial-time Turing machine\}.
\end{defn}
In particular, when $f \in {\sf fP}$ then ${\sf Dom}(f)$ is in {\sf P}. 
The set {\sf fP} is closed under composition of functions, and the identity 
function is in {\sf fP}, so {\sf fP} is a monoid. 
The above characterization of {\sf P} $\neq$ {\sf NP} now becomes: 

\smallskip

 \ \ \ \ \ \ {\sf P} $\neq$ {\sf NP} \, 
{\it iff \, {\sf fP} is a non-regular monoid.}

\medskip

We will use some more definitions that apply to functions. 

For any function $f: X \to Y$, the equivalence relation ${\sf mod}_f$ on
${\sf Dom}(f)$ is defined by: $\, x_1 \ {\sf mod}_f \ x_2 \, $ iff
$\, f(x_1) = f(x_2)$.
The set of equivalence classes of $\, {\sf mod}_f \,$ is
$\, \{f^{-1}(y) : y \in {\sf Im}(f)\}$.

A {\em choice set} for $f$ is, by definition, a subset of 
${\sf Dom}(f)$ that contains exactly one element of each 
${\sf mod}_f$-class. 
A {\em choice function} for $f$ is, by definition, any inverse $f'$ of $f$ 
such that ${\sf Dom}(f') = {\sf Im}(f)$. A choice function $f'$ maps each 
element $y \in {\sf Im}(f)$ injectively to an element of $f^{-1}(y)$, and 
every ${\sf mod}_f$-class contains exactly one element of ${\sf Im}(f')$; 
so ${\sf Im}(f')$ is a choice set. 
For a choice function $f'$ we have 
$\, f \circ f' = f \circ f'|_{{\sf Im}(f)} = {\sf id}_{{\sf Im}(f)}$, hence
$\, f f' f = f$; we also have $f' f f' = f'$; more generally, if 
$f'_1, f'_2$ are two choice functions for the same function $f$ then 
$\, f'_2 f f'_1 = f'_2$.
A choice function for $f$ is uniquely determined by $f$ and a choice set. 

A {\em representative choice function} for $f$ is, by definition, any 
function $r$ such that ${\sf Im}(r) \subseteq {\sf Dom}(r) = {\sf Dom}(f)$, 
${\sf Im}(r)$ is a choice set for $f$, ${\sf mod}_r = {\sf mod}_f$, and 
$r|_{{\sf Im}(r)} = {\sf id}_{{\sf Im}(r)}$.  Thus, $r$ maps each 
${\sf mod}_f$-class $[x]_f$ to one and the same chosen element of $[x]_f$. 
It follows that $r$ is an idempotent; and if $r_1, r_2$ are two 
representative choice functions for the same function $f$ then 
$\, r_2 \circ r_1(.) = r_2(.)$.
If $f \in {\sf fP}$ and if in addition $r \in {\sf fP}$, then 
$f \equiv_{\cal L} r$ (by Prop.\ 2.1 in \cite{s1f}).
So in that case $f$ is a regular element of {\sf fP}, since it is 
$\cal L$-equivalent to an idempotent. 
Hence, if every function in {\sf fP} had a representative choice function 
in {\sf fP}, then {\sf fP} would be regular, which would imply 
${\sf P} = {\sf NP}$.

There is a one-to-one correspondence between the choice functions of $f$ 
and the representative choice functions:
If $c$ is a choice function then $\rho_c(.) = c \circ f(.)$ is a 
representative choice function; and if $r$ is a representative choice 
function of $f$ then  $v_r(.) = r \circ f^{-1}(.)$ is a choice 
function (where $f^{-1}$ is the inverse relation of $f$, and $\circ$ is 
composition of relations).
Moreover, if $r$ is a representative choice function of $f$ then
$\rho_{v_r} = r$; and if $c$ is a choice function of $f$ then 
$v_{\rho_c} = c$.  

\medskip

For functions in general, the existence of a choice function, 
the existence of a representative choice function, and 
the existence of a choice set,
are equivalent to the axiom of choice.
Because of this connection we are especially interested in inverses $f'$ 
of $f$ that are choice functions, i.e., that satisfy 
${\sf Dom}(f') = {\sf Im}(f)$.

These three formulations of the axiom of choice can also be considered for 
{\sf fP}: 
(1) ``Every $f \in {\sf fP}$ has an inverse in {\sf fP}.'' 
Equivalently, ``$f$ has a choice function in {\sf fP}'', and also 
equivalently, ``$f$ has a mutual inverse in {\sf fP}''.
(2) ``Every $f \in {\sf fP}$ has a representative choice function in 
{\sf fP}.''
(3) ``Every $f \in {\sf fP}$ has a choice set in {\sf P}.''
 \, It is an open problem problem whether the three statements are true. 

We saw that (1) is true iff ${\sf P} = {\sf NP}$. Moreover, (1) implies (2); 
indeed, if $f$ has an inverse $f' \in {\sf fP}$ then $f'f$ is a 
representative choice function in {\sf fP}.
And (2) implies (3); indeed, if $r \in {\sf fP}$ is a representative choice
function of $f$, then $r$ is regular (being an idempotent), hence 
${\sf Im}(r)$ belongs to {\sf P} by Prop.\ 1.9 in \cite{s1f}; and 
${\sf Im}(r)$ is a choice set for $r$ and for $f$. It remains an open 
problem whether other implications between (1), (2), (3) hold. 
Nevertheless, we have the following.

\begin{pro} \label{existsInj1wf}
 \ There exists an {\em injective} one-way function (for worst-case 
complexity) iff there exists a one-way function that has a choice set in 
{\sf P}.
\end{pro}
{\bf Proof.} If $g \in {\sf fP}$ is injective then ${\sf Dom}(g)$ is a
choice set for $g$, and ${\sf Dom}(g) \in {\sf P}$. If $g$ is also one-way,
then it is therefore a one-way function with a choice set in {\sf P}. 

The converse follows immediately from the next Lemma. 
 \ \ \ $\Box$

\begin{lem} \label{Inj1wfVSchoice}
 \ If $f \in {\sf fP}$ has a choice set $C \in {\sf P}$, but $f$ is not
regular (i.e., (3) holds but (1) does not), then the restriction $f|_C$ is 
an injective one-way function.
\end{lem}
{\bf Proof.} Since $f \in {\sf fP}$ and $C \in {\sf P}$ then we have 
$f|_C \in {\sf fP}$. Moreover, if $C$ is a choice set for $f$ then $f|_C$ 
is injective. And if $C$ is a choice set, any inverse of $f|_C$ is also an 
inverse of $f$. Hence $f|_C$ has no inverse in {\sf fP} (since $f$ is not 
regular).
 \ \ \ $\Box$

\bigskip

The motivation for this paper is based on the following simple observation:
If $f: A^* \to A^*$ is any partial function and if $f'$ is an inverse of
$f$, then the restriction $f'|_{{\sf Im}(f)}$ is an {\em injective} 
inverse of $f$. 
Moreover, if $f \in {\sf fP}$ and $f' \in {\sf fP}$, then
${\sf Im}(f) \in {\sf P}$ (by Prop.\ 1.9 in \cite{s1f}); hence
$f'|_{{\sf Im}(f)} \in {\sf fP}$.
Thus we have: 
$\, {\sf P} \neq {\sf NP} \, $ iff there exists 
$f \in {\sf fP}$ such that $f$ has no {\em injective} inverse in {\sf fP}.

The classes {\sf P} and {\sf NP} are defined in terms of {\em sets} of 
strings (``formal languages'').  To add more structure we use 
{\em functions}, and we characterized {\sf P} and {\sf NP} by properties of  
{\em monoids} of functions \cite{s1f}.  Our next step, in the present 
paper, is to characterize {\sf P} versus {\sf NP} by properties of 
{\em inverse monoids}, and of {\em groups}.
However, whether this approach will help solve the {\sf P} versus {\sf NP} 
problem, remains to be seen.

\medskip

\bigskip

\noindent  {\bf Overview:}
  In Section {\bf 2} we introduce the monoid {\sf invfP} consisting of all 
injective regular functions in {\sf fP}. 
We prove that every regular function in {\sf fP} has an inverse in 
{\sf invfP}, and that {\sf invfP} is a maximal inverse submonoid of 
{\sf fP}.

  In Section {\bf 3} we show that the polynomial-time polynomially balanced 
{\em injective} Turing machines form a machine model for {\sf invfP}, i.e., 
that a function $f$ belongs to {\sf invfP} iff $f$ is computed by such a 
Turing machine. We conclude from this that {\sf invfP} is finitely generated. 
We also consider polynomial-time polynomially balanced injective Turing 
machines with an {\sf NP}-oracle, and we show that the set 
${\sf invfP}^{\sf (NP)}$ of functions computed by such Turing machines is a 
finitely generated inverse monoid. To prove the latter, we show that there 
exist languages that are {\sf NP}-complete with respect to one-one 
reductions in {\sf invfP}.  We show that every function in {\sf fP} has 
an inverse in ${\sf invfP}^{\sf (NP)}$.

  In Section {\bf 4} we show that ${\sf invfP} \neq {\sf invfP}^{\sf (NP)}$ 
iff ${\sf P} \neq {\sf NP}$. 
We introduce the set {\sf cofP} of all functions in 
${\sf invfP}^{\sf (NP)}$ that have an inverse in {\sf fP}, i.e., that are a
co-inverse of a function in {\sf fP}.
We show that this is a finitely generated monoid.  We prove that 
${\sf P} \neq {\sf NP}$ iff ${\sf invfP} \neq {\sf cofP}$, iff 
${\sf cofP} \neq {\sf invfP}^{\sf (NP)}$, iff {\sf cofP} is not 
regular. 

We also introduce the monoid {\sf injfP} of all injective functions in 
{\sf fP}; $\, {\sf injfP}$ is also equal to 
${\sf invfP}^{\sf (NP)} \, \cap \, {\sf fP}$.  We show that the inverse 
monoid {\sf invfP} is the set of the regular elements of {\sf injfP}.  
Hence, injective one-way functions (for worst-case complexity) exist iff 
${\sf invfP} \neq {\sf injfP}$.
We do not know whether {\sf injfP} is finitely generated; if it turns out 
that {\sf injfP} is not finitely generated then ${\sf P} \neq {\sf NP}$. 

  In Section {\bf 5} we show that every element of {\sf fP} is equivalent 
(with respect to inversive-reductions) to an element of {\sf fP} that has 
an inverse in a {\em subgroup} of ${\sf invfP}^{\sf (NP)}$. 

\bigskip

\bigskip

\noindent  {\bf Notation for the monoids used, and their definition:} 

\medskip

\noindent \begin{tabular}{lll}
{\sf fP}    &   Definition \ref{DEFfP} & (polynomially balanced 
 polynomial-time computable functions)  \\   
{\sf invfP} & Definition \ref{DEFinvfP} & (injective functions in {\sf fP}
 that have an inverse in {\sf fP}) \\  
${\sf invfP}^{\sf (NP)}$ & Definition \ref{DEFinvfPNP} & (functions
 computed by injective Turing machines with {\sf NP}-oracle) \\   
{\sf cofP} \ \ \ \  &  Definition \ref{DEFcofP} & (functions in 
 ${\sf invfP}^{\sf (NP)}$ that have an inverse in {\sf fP}) \\  
{\sf injfP} & Definition \ref{DEFinjfP} & (injective functions in {\sf fP})
\end{tabular} 

\medskip

\noindent \ The relation between these monoids is shown in Figure 1 
(after Prop.\ \ref{invfPNPcapfP}).

%%%%%%%%%%%%%%%%%%%%%%%%%%%%%%%%%%%%%%%%%%%%%%%%%%%%%%%%%%%%%%%%%%
\section{Injective inverses of regular elements}

For any function $f: A^* \to A^*$, we will be interested in inverses $f'$ of 
$f$ that satisfy ${\sf Dom}(f') = {\sf Im}(f)$; such inverses are 
necessarily injective. Moreover, when ${\sf Dom}(f') = {\sf Im}(f)$ we not 
only have $f f' f = f$, but also $f' f f' = f'$. So $f$ and $f'$ are mutual 
inverses. 

An injective inverse $f'$, as above, has exactly one inverse $f''$ that
satisfies ${\sf Dom}(f'') = {\sf Im}(f')$, namely $f'' = f'^{-1}$ (the 
set-theoretic inverse function of $f'$); moreover, 
${\sf Im}(f'^{-1}) = {\sf Im}(f)$.  
And $f'^{-1} = f|_{{\sf Im}(f')}$, i.e., the restriction of $f$ to the
choice set ${\sf Im}(f')$ of $f'$.

Here are some more simple facts about inverses:
 \ If $f'$ is a mutual inverse of $f$ then
$\, {\sf Im}(f') \subseteq {\sf Dom}(f)$. 
Moreover, $f'$ is a mutual inverse of $f \, $ iff 
$\, {\sf Im}(f') = {\sf Im}(f'|_{{\sf Im}(f)})$.
Note that $f'|_{{\sf Im}(f)}$ is the choice function determined
by $f'$, so ${\sf Dom}(f'|_{{\sf Im}(f)}) = {\sf Im}(f)$. Thus 
an inverse $f'$ is a mutual inverse of $f$ \ iff \ ${\sf Im}(f')$ is
the choice set determined by $f'$ in ${\sf Dom}(f)$.
As a consequence, a mutual inverse $f'$ of $f$ is injective iff
${\sf Dom}(f') = {\sf Im}(f)$.

\begin{pro} \label{prod_inverses}
 \ If $g_1$ and $g_2$ are injective regular elements of {\sf fP}, then
$g_2 \circ g_1$ is also injective and regular in {\sf fP}. 
\end{pro}
{\bf Proof.} The composite of injective functions is obviously injective. 
Since $g_i$ is regular (for $i = 1, 2$), it has the injective regular
function $g_i^{-1}$ as an inverse. Indeed, ${\sf Im}(g_i) \in {\sf P}$, 
and $g_i^{-1} = g_i'|_{{\sf Im}(g_i)}$ for any inverse $g_i' \in {\sf fP}$
of $g_i$; hence, $g_i^{-1} \in {\sf fP}$.  
Thus, the injective function $g_1^{-1} \circ g_2^{-1}$ is a mutual inverse 
of $g_2 \circ g_1$, and 
${\sf Dom}(g_1^{-1} \circ g_2^{-1}) = {\sf Im}(g_2 \circ g_1)$. 
 \ \ \ $\Box$

\bigskip

\noindent Note that in general, the product of regular elements in {\sf fP}
need not be regular (unless ${\sf P} = {\sf NP}$); but by the above 
proposition, the product of {\em injective} regular elements is regular.

\begin{defn} \label{DEFinvfP}
 \ Let {\sf invfP} denote the set of injective regular elements of {\sf fP}.
\end{defn}

\begin{cor} \label{injreg}
 \ The set {\sf invfP} of injective regular elements of {\sf fP} is an 
{\em inverse monoid}. Every regular element of {\sf fP} has an inverse in 
{\sf invfP}.  \ \ \ \ \ $\Box$
\end{cor}
Corollary \ref{injreg} implies that ${\sf P} = {\sf NP}$ iff
every element of {\sf fP} has an inverse in {\sf invfP}. 

Another motivation for {\sf invfP} will be seen in Prop.\ 
\ref{NPcomplINVred}, where is is shown that there exist languages that are 
{\sf NP}-complete with respect to one-one reductions in {\sf invfP}.

\medskip

A monoid $M_1$ is called an {\em inverse submonoid} of a monoid $M$ iff 
$M_1$ is submonoid of $M$, and $M_1$ is an inverse monoid by itself.
So, an element of $M_1$ has exactly one mutual inverse in $M_1$, but it 
could have additional mutual inverses in $M$.

\begin{pro} \label{invfP_max}
 \ The inverse submonoid ${\sf invfP}$ is a {\em maximal} inverse submonoid
of {\sf fP}; i.e., if $M$ is an inverse monoid such that
$\, {\sf invfP} \subseteq M \subseteq {\sf fP}$, then $\, {\sf invfP} = M$.
\end{pro}
{\bf Proof.} Let $M$ be an inverse submonoid such that 
${\sf invfP} \subseteq M \subseteq {\sf fP}$.  For a contradiction, let us 
assume that ${\sf invfP} \neq M$. Then $M$ contains a non-injective element
$f$, which is regular (since $M$ is inverse). Since $f$ is regular, $f$ has 
a choice function $f'_1$; since $f$ is non-injective, $f$ has at least one 
other choice function $f'_2$. 

Moreover, since $f'_1 \in {\sf invfP}$, we can pick $f'_2$ in such a way
that $f'_2 \in {\sf invfP}$ too.
Indeed, let $x_1, x_2 \in {\sf Dom}(f)$ be such that $f(x_1) = f(x_2)$, and 
$x_1$ is in the choice set ${\sf Im}(f_1')$, and 
$x_2 \not\in {\sf Im}(f_1')$. Let $\tau$ be the transposition of $x_1$ and
$x_2$ (and $\tau$ is the identity elsewhere). Then $f_2' = f_1' \circ \tau$
belongs to {\sf invfP}, and is the same as $f_1'$, except that $x_2$ has 
replaced $x_1$ in the choice set. 

Since the choice functions $f_1', f_2'$ belong to {\sf invfP} 
($\subseteq M$), $f$ has two mutual inverses in $M$, which implies that 
$M$ is not an inverse monoid. 
 \ \ \ $\Box$

\medskip

\noindent {\bf Remarks.}  (1) Prop.\ \ref{invfP_max} holds, whether {\sf fP} 
is regular or not.  \\      
(2) An inverse submonoid of {\sf fP} need not consist of injective functions
only, and the submonoid {\sf invfP} is not the only maximal inverse 
submonoid of {\sf fP}. E.g., {\sf fP} contains non-injective idempotents, 
and these are contained in maximal inverse monoids that are different from 
{\sf invfP} (since they contain non-injective elements).

\begin{pro} \label{maxinvDisjoint} 
 \ A maximal subgroup of {\sf fP}, and more generally, any $\cal L$-class of
{\sf fP}, is either {\em disjoint} from {\sf invfP} or entirely contained in 
{\sf invfP}.
\end{pro}
{\bf Proof.} An $\cal L$-class of {\sf fP} that intersects {\sf invfP} 
contains injective elements, hence it consist entirely of injective elements
(by Prop.\ 2.1 in \cite{s1f}). Also, an $\cal L$-class of {\sf fP} that
intersects {\sf invfP} contains regular elements, hence it consist entirely 
of regular elements (it is a well-known fact from semigroup theory that if
a $\cal D$-class contains regular elements then it consists entirely of
regular elements; see e.g.\ \cite{CliffPres, Grillet}).
Hence, this $\cal L$-class consists entirely of regular injective elements,
hence it is contained in {\sf invfP}.
 \ \ \ $\Box$

\medskip

On the other hand, every regular $\cal R$-class of {\sf fP} intersects both
{\sf invfP} and ${\sf fP} - {\sf invfP}$.  Indeed, for every regular element 
$f \in {\sf fP}$ we have $f \equiv_{\cal R} {\sf id}_{{\sf Im}(f)} \,$ (and
${\sf Im}(f)$ is in {\sf P} when $f$ is regular).  But a regular 
$\cal R$-class always contains some non-injective elements (by Prop.\ 2.1 in 
\cite{s1f}).

%%%%%%%%%%%%%%%%%%%%%%%%%%%%%%%%%%%%%%%%%%%%%%%%%%%
\section{Computing inverses}

There is a machine model that exactly characterizes the injective inverses 
of the regular elements of {\sf fP}, namely the polynomially balanced 
polynomial-time injective Turing machines. 
As a consequence we will prove that {\sf invfP} is finitely generated. 

Turing machines of the above type, with an {\sf NP}-oracle added, form a
model of computation for another inverse monoid, ${\sf invfP}^{\sf (NP)}$ 
(given in Def.\ \ref{DEFinvfPNP}), which contains some injective inverses 
for every element of {\sf fP}.

%%%%%%%%%%%
\subsection{Injective Turing machines}

There is a very simple machine model for the elements of {\sf invfP}, namely
the polynomially balanced polynomial-time injective Turing machines. 
A deterministic Turing machine is called {\em injective} iff the transition 
table of the Turing machine describes a (finite) injective function; for 
details, see \cite{Bennett73, Bennett89}.  

The {\em reverse} of a Turing machine $M$ is the machine obtained by 
reversing every transition of $M$ (and also switching start and accept
states). The reverse of a deterministic Turing machine $M$ is is not 
deterministic, unless $M$ is injective.  

\smallskip

\noindent {\bf Remark.} In the literature, injective Turing machines are 
called ``reversible'', because of historic connections with the study 
of computation as a ``reversible process'' (in the sense of thermodynamics).
However, calling an injective Turing machine ``reversible'' can be 
misleading, since the transitions obtained by reversing an injective Turing 
machine $M$ are not part of $M$ (but of a different machine).

\begin{pro} \label{noninvertibleInj} \hspace{-.13in} {\bf .}

\noindent {\bf (1)} \ The reverse of a polynomially balanced 
polynomial-time injective Turing machine is also a polynomially balanced 
polynomial-time injective Turing machine. 
If a function $f$ is polynomially balanced and it is computed by a 
polynomial-time injective Turing machine, then $f^{-1}$ is
also polynomially balanced and it is computed by a polynomial-time
injective Turing machine.     \\   
{\bf (2)} \ There exists a polynomial-time injective Turing machine (that is
not polynomially balanced)  whose reverse does not have polynomial 
time-complexity. 
There exists an injective function $g$ that is computed by a polynomial-time 
injective Turing machine, whose inverse $g^{-1}$ is not computable by a 
polynomial-time Turing machine.  
\end{pro}
{\bf Proof.} (1) An injective Turing machine $M$, computing an injective 
function $f$, can be run in reverse. This yields a new injective Turing
machine $M'$, computing $f^{-1}$. If $M$ runs in polynomial time, with 
polynomial $p_T(.)$, and it is polynomially balanced, with polynomial 
$p_B(.)$, then $M'$ is also polynomially balanced, with polynomial 
$p_B(.)$ too, and runs in polynomial time, with polynomial bound 
$p_T \circ p_B(.)$.

\smallskip

\noindent (2) An example of such a function is

\smallskip
  
 \ \ \  \ \ \  $g: \ a^{(2^m)} \ \longmapsto \ a^m$, \ \ for all 
                                            $m \in {\mathbb N}$,

\smallskip

\noindent where $a$ is a fixed letter; $\, g$ is undefined on any input 
that is not of the form $a^k$ with $k$ a power of 2.

Obviously, $g$ is injective and polynomial-time computable, but not 
polynomially balanced. Hence, $g^{-1}$ is not computable in polynomial 
time, since its output is exponentially longer than its input. 
However, $g$ is computable in linear time by the following injective Turing
machine: We use a Turing machine with a {\em rubber tape} for the input, 
and an ordinary tape for the output. A rubber tape is a tape on which one 
can not only replace one letter by another one (as on an ordinary tape), but
where one also can insert a letter or remove a letter (in one transition).
A rubber tape can easily be simulated by two stacks; a stack is a special
case of an ordinary tape.
The machine that computes $g$ works as follows, in outline.
It has a main loop (``{\tt while} $\ldots$''); it has an inner loop that is
executed at the beginning of the body of the main loop. For this inner loop,
the machine uses two states to count the input positions modulo 2, and 
erases every second $a$ (which is possible on a rubber tape). When the right 
end of the input tape is reached, the state must correspond to an even 
number of letters $a$ on the input tape (otherwise, the machine rejects and 
has no output for this input). 
At this moment, one $a$ is printed on the output tape.

\medskip

\hspace{.1in} {\tt while} the input tape is not empty: \verb|  // main loop|

\hspace{.1in} \{ \hspace{.23in} 
               in a loop, erase every second $a$ on the input tape;
                            \verb|  // inner loop|

\hspace{.5in}       when the right end of the input tape is reached,
                            \verb|  // after the inner loop|

\hspace{0.5in}          {\tt if} the number of letters $a$ read in the 
                                inner loop was even:  

\hspace{.8in}          {\tt then} add an $a$ on the output tape;

\hspace{.8in}          {\tt else} reject;

\hspace{.5in}       in a loop, move the head of the input tape back to 
                        the left end; \verb|  // 2nd inner loop|

\hspace{.1in} \}

\medskip
 
\noindent This program runs in linear time, and every step is injective.

This injective Turing machine for $g$ can be reversed, which yields an 
injective Turing machine for $g^{-1}$ with exponential time-complexity.   
 \ \ \ $\Box$

\begin{pro} \label{invfPinjTM}
 \ Let $h$: $A^* \to A^*$ be an injective function.  Then 
$h$ belongs to {\sf invfP} \ iff \ $h$ is polynomially balanced, and is 
computable by some polynomial-time injective Turing machine.
\end{pro}
{\bf Proof.} 
For any $h \in {\sf invfP}$, $h$ is injective and computable in
polynomial time; moreover, $h^{-1} \in {\sf invfP}$, i.e., $h^{-1}$ is also
injective and computable in polynomial time. Hence, Bennett's theorem is 
applicable (see \cite{Bennett73, Bennett89}), so $h$ can be computed by a 
polynomial-time {\em injective} deterministic Turing machine. Also, since
both $h$ and $h^{-1}$ are computable in polynomial time, $h$ is polynomially
balanced. 
Conversely, if $h$ is polynomially balanced and is computed by a 
polynomial-time injective Turing machine, then $h^{-1}$ can also be computed 
by such a machine, by Prop.\ \ref{noninvertibleInj}.
 \ \ \ $\Box$

\bigskip

We saw in the Introduction that $f \in {\sf fP}$ is regular iff 
$f$ has an inverse in {\sf invfP}. Hence by Prop.\ \ref{invfPinjTM}, we 
have for all $f \in {\sf fP}$:

{\it $f$ is regular \ iff \ $f$ has an inverse that is computable by a 
polynomially balanced polynomial-time injective Turing machine. 
And ${\sf P} = {\sf NP}$ \ iff \ every $f \in {\sf fP}$ has an inverse that 
can be computed by a polynomially balanced polynomial-time injective Turing
machine.}

%%%%%%%%%%%
\subsection{Evaluation functions and finite generation}

Based on the transition table of a Turing machine, one can easily check
whether this Turing is deterministic, and whether it is injective.
Moreover, we saw in \cite{s1f} and \cite{infgen} that one can add a built-in
polynomial time bound and balance bound into the transition table.  
This can be done for injective Turing machines too, without destroying
injectiveness; indeed, a time bound $p(|x|)$ can be computed injectively on 
input $x$ (where $p$ is a stored polynomial, described by its degree and 
coefficients).  
So we can design a set of strings (programs) that describe all polynomially 
balanced polynomial-time injective Turing machines. And just as 
in \cite{s1f}, for every polynomial $q(n) = a \cdot (n^k + 1)$, there exists
an evaluation function ${\sf ev}_q$ that evaluates all injective Turing 
machine programs with built-in polynomial less than $q$.   
The details are the same as in \cite{s1f}, and injectiveness doesn't change
any reasoning.
We call such a Turing machine description an {\sf invfP}-{\it program}.

A Turing machine with program $w$ will be denoted by $M_w$; we denote the 
injective input-output function of $M_w$ by $\phi_w$. 

An evaluation function maps a pair $(w,x)$, consisting of a program $w$ and 
an input $x$ (for $M_w$), to $(w, \phi_w(x))$, i.e., to that same program 
and the program-output. In order to represent a pair of strings by one 
string we use the prefix code $\{00, 01, 11\}$, and the function 
{\sf code}(.) defined by 

\smallskip

 \ \ \   ${\sf code}(0) = 00$, \ \ ${\sf code}(1) = 01$. 

\smallskip

\noindent The pair of strings $(w,x)$ is represented unambiguously by 
${\sf code}(w) \, 11 \, x$, where 11 acts as a separator (since 
$\{00, 01, 11\}$ is a prefix code).
The {\em evaluation function} ${\sf injEv}_q$ is defined by 

\smallskip

 \ \ \  ${\sf injEv}_q\big({\sf code}(w) \, 11 \, x \big) \ = \ $
${\sf code}(w) \, 11 \ \phi_w(x)$

\smallskip

\noindent for any {\sf invfP}-program $w$ with built-in polynomial less than
$q$, and any $x \in {\sf Dom}(\phi_w)$.
The function ${\sf injEv}_q$ is itself injective and polynomially balanced 
and polynomial-time computable (with a larger polynomial than $q$ however).
Moreover, ${\sf injEv}_q$ is regular; indeed, the unique mutual inverse of 
${\sf injEv}_q$ is ${\sf injEv}_q^{-1}$, and this belongs to {\sf fP} since
$\phi_w(x)$ is injective and regular, i.e., $\phi_w^{-1} \in {\sf fP}$ for 
every {\sf invfP}-program $w$.
So, ${\sf injEv}_q \in {\sf invfP}$.

\medskip

The proof of Prop.\ 4.5 in \cite{s1f}, showing that {\sf fP} is finitely 
generated, goes through without much change. We use the relation 

\medskip

\noindent $(\star)$
\hspace{.5in} $\phi_w(x) \ = \ $
$\pi_{_{2 \, |w'| + 2}}' \circ {\sf contr} \circ {\sf recontr}^{2 \, m}$
$\circ$ ${\sf injEv}_q$ $\circ$
${\sf reexpand}^{2m} \circ {\sf expand}$  $\circ$
$\pi_{_{{\sf code}(w) \, 11}}(x)$,

\medskip

\noindent where $w$ is a program of a polynomially balanced 
polynomial-time injective Turing machine with built-in polynomial $< q$, 
and $\, w' = {\sf co}^{2m+1} \circ {\sf ex}^{2m+1}(w)$.
See Section 4 of \cite{s1f} for the definition of 
{\sf contr}, {\sf recontr}, {\sf reexpand}, {\sf expand}, 
$\pi_n'$, $\pi_v$, {\sf co}, and {\sf ex}. The role of {\sf expand} and 
{\sf reexpand} is to increase the input length, so as to reduce complexity 
by a padding argument; when complexity is below $q$, ${\sf injEv}_q$ can be
applied; after that, {\sf recontr} and {\sf contr} remove the padding.

The functions {\sf contr}, {\sf recontr}, ${\sf injEv}_q$, {\sf reexpand},
{\sf expand}, and $\pi_{_{{\sf code}(w) \, 11}}$ are injective and regular
(i.e., they belong to {\sf invfP}).
In the relation $(\star)$ we can replace $\pi_{_{2 \, |w'| + 2}}'$ by the
injective regular function 
$\, \pi_{_{{\sf code}(w) \, 11}}': \ {\sf code}(w) \, 11 \ x \ $
$\longmapsto \ x$. \ Moreover, $\pi_{_{{\sf code}(w)}}$ is generated by 
$\{\pi_0, \pi_1\}$, and $\pi_{_{2 \, |w'| + 2}}'$ is generated by
$\{\pi_0', \pi_1'\}$.

Thus, $\phi_w$ is generated by the finite set of functions 
$\, \{{\sf contr},$ ${\sf recontr},$ ${\sf injEv}_q,$ ${\sf reexpand},$
${\sf expand},$ $\pi_0,$ $\pi_1,$ $\pi_0',$ $\pi_1'\}$ $\subset$ 
${\sf invfP}$, and we have:

\begin{pro} \label{invfP_FinGen}
 \ The inverse monoid {\sf invfP} is finitely generated, as a monoid.
 \ \ \  \ \ \ $\Box$
\end{pro}

For every {\sf invfP}-program $w$, we can easily obtain an 
{\sf invfP}-program (let's call it $w'$) for $\phi_w^{-1}$; for this 
purpose we simply reverse the injective Turing machine described by $w$. 
Thus the function $w \mapsto w'$ is polynomially balanced and 
polynomial-time computable; 
moreover, $w \mapsto w'$ is injective and involutive, so the function
$w \mapsto w'$ is in {\sf invfP}. Thus we proved:

\begin{pro} \label{inv_mapforInvfP}
 \ There exists a function ${\sf prog}_{\sf inv} \in {\sf invfP}$ such 
that ${\sf prog}_{\sf inv}(w)$ is an {\sf invfP}-program for $\phi_w^{-1}$,
for every {\sf invfP}-program $w$.
  \ \ \  \ \ \ $\Box$
\end{pro}

We would like to extend the above function ${\sf prog}_{\sf inv}$ to all 
``regular'' {\sf fP}-programs, i.e., the {\sf fP}-programs $w$ for which 
$\phi_w$ is regular in {\sf fP}. We can do this by using 
{\em universal search}, a.k.a.\ {\em Levin search}; this is described in
\cite{Levin73, Levin84}, but without much detail; for a detailed exposition, 
see for example \cite{LiVit} (Theorem 7.21 in the 1993 edition).
In the general regular case, ${\sf prog}_{\sf inv}(w)$ is just a program 
such that $\phi_{{\sf prog}_{\sf inv}(w)}$ has polynomial 
time-complexity, but ${\sf prog}_{\sf inv}(w)$
does not have a built-in polynomial for its time-complexity; Levin
search is not able to explicitly find such a polynomial, although it exists
if $\phi_w$ is regular. 

\begin{pro} \label{Levinsearch} {\bf (inversion by Levin search).}
 \ There exists a function ${\sf prog}_{\sf inv} \in {\sf fP}$ such that for 
every {\sf fP}-program $w$, ${\sf prog}_{\sf inv}(w)$ is a program for a 
mutual inverse of $\phi_w$ satisfying 
${\sf Dom}(\phi_{{\sf prog}_{\sf inv}(w)})$ $=$ ${\sf Im}(\phi_w)$. 

\smallskip

The time-complexity of $\, \phi_{{\sf prog}_{\sf inv}(w)} \,$ is
$\, \Theta(T_{\phi_w}) + \Theta(T_w')$, where $\Theta(T_{\phi_w})$ is (up 
to big-$\Theta$) the optimal time-complexity of all Turing machines for 
$\phi_w$, and $\Theta(T_w')$ is (up  to big-$\Theta$) the optimal 
time-complexity of all Turing machines for all mutual inverses of $\phi_w$.

\smallskip

\noindent In particular: 

If $\phi_w$ is regular then $\phi_{{\sf prog}_{\sf inv}(w)}$ has polynomial 
time-complexity. But the program ${\sf prog}_{\sf inv}(w)$ does not have a 
built-in polynomial for its time-complexity, i.e., ${\sf prog}_{\sf inv}(w)$ 
is not an {\sf fP}-program.

If $\phi_w$ is not regular, ${\sf prog}_{\sf inv}(w)$ is a 
non-polynomial-time program.  The time-complexity of 
$\phi_{{\sf prog}_{\sf inv}(w)}$ has nevertheless an exponential upper bound.
\end{pro}
{\bf Proof.} In our version of universal search the input has the form
$(w,y)$, where $w$ is an {\sf fP}-program, and $y$ is a possible output of 
the Turing machine $M_w$. 
Remark: For convenience we write $(w,y)$ as a pair of words, but our 
universal search will actually use the single word 
$\, {\sf code}(w) \, 11 \, y$.

The output of universal search is $(w,x)$ such that 
$x \in \phi_w^{-1}(y)$, if $w$ is an {\sf fP}-program, and 
$y \in {\sf Im}(\phi_w)$ (i.e., if $\phi_w^{-1}(y) \neq \varnothing$); 
there is no output if $w$ is not an {\sf fP}-program, or if 
$y \not\in$ ${\sf Im}(\phi_w)$. 
Thus universal search computes a mutual inverse ${\sf ev}'$ of the general 
evaluation function ${\sf ev}$ for {\sf fP}. Both ${\sf ev}'$ and ${\sf ev}$ 
are partial recursive (but they do not belong to {\sf fP} since no 
polynomial bound is prescribed, as for ${\sf ev}_q$). 
Since we restrict universal search to {\sf fP}-programs (hence, with 
built-in polynomial complexity bound), the domain of ${\sf ev}$ is decidable;
so universal search is an algorithm (that always halts).
Indeed, by rejecting programs with no polynomial time bound, the search 
algorithm also rejects programs that do not halt.  And for a polynomial
program with known polynomial complexity, inversion has a predictable 
exponential-time upper-bound.  

For a fixed $w$ in the input, universal search computes
${\sf ev}'(w, \cdot) = (w, \phi_w'(.))$, where $\phi_w'(.)$ is a mutual 
inverse of $\phi_w$. For the details of the universal search algorithm we
refer to \cite{LiVit} (proof of Theorem 7.21, p.\ 412 in the 1993 edition).
From the specification of universal search and $w$, we immediately derive a 
program for ${\sf ev}'(w, \cdot)$; this program is called 
${\sf prog}_{\sf inv}(w)$. Since ${\sf prog}_{\sf inv}(w)$ is easily 
obtained from $w$, the function ${\sf prog}_{\sf inv}(.)$ belongs to 
{\sf fP}.

The program ${\sf prog}_{\sf inv}(w)$ has minimum time-complexity
$\, \Theta(T_{\phi_w}) + \Theta(T_w')$; this is proved in \cite{LiVit}.  
It follows that the time-complexity of ${\sf prog}_{\sf inv}(w)$ has an 
exponential upper bound (since $w$ is an {\sf fP}-program).  
It also follows that ${\sf prog}_{\sf inv}(w)$ is a polynomial-time 
program if $\phi_w$ is regular, i.e., if $\phi_w$ has a polynomial-time 
inverse.
But universal search does not find out explicitly what that polynomial is, 
so ${\sf prog}_{\sf inv}(w)$ cannot have a built-in polynomial for its 
time-complexity bound; i.e., ${\sf prog}_{\sf inv}(w)$ is not an 
{\sf fP}-program. 
 \ \ \ $\Box$

%%%%%%%%%%%%
\subsection{Injective Turing machines with {\sf NP}-oracle}

In preparation for the next Section, we generalize polynomially balanced 
polynomial-time injective Turing machines by adding an oracle from {\sf NP}.

\begin{defn} \label{DEFinvfPNP} 
 \ By $\, {\sf invfP}^{\sf (NP)}$ we denote the set of all functions computed 
by polynomially balanced polynomial-time {\em injective} Turing machines, 
with an oracle belonging to {\sf NP}. 
\end{defn}
By Prop.\ \ref{invfPinjTM} we have:
 \ ${\sf invfP} \subseteq {\sf invfP}^{\sf (NP)}$.

\bigskip

An injective Turing machine with an oracle is a deterministic and injective
machine. Indeed, it computes deterministically and injectively inbetween 
oracle transitions. And in an oracle transition the only change is in the 
state, which goes from the query state $q_{\sf qu}$ to either $q_{\sf yes}$ 
or $q_{\sf no}$; the contents and head positions of the tapes (including the 
query tape) remain unchanged. 
The next transition will go to a different state than $q_{\sf yes}$ or 
$q_{\sf no}$, so the state $q_{\sf yes}$ (or $q_{\sf no}$) uniquely 
determines the previous state $q_{\sf qu}$. Thus, an oracle transition is 
deterministic and injective. 

\smallskip

The reverse of an injective Turing machine with an oracle has a slightly 
different format than an an injective Turing machine with an oracle, as
described above. Indeed, in an oracle call the following transition happens:
The machine is in the query state $q_{\sf qu}$ and then it enters an 
answer state, either $q_{\sf yes}$ or $q_{\sf no}$, according as the word on
the query tape belongs to the oracle set or not. 
(Before a query, a word is written on the query tape; this happens in the
course of a possibly long computation. 
After a query, the computation continues, and during this computation the 
query tape content may be gradually erased or changed.) 

When this sequence is reversed, an answer state occurs before the query
state. From the answer state, the reverse computation goes to the query 
state, provided that the answer state is $q_{\sf yes}$ and the query word is 
in the oracle language, or if the answer state is $q_{\sf no}$ and the query
word is not in the oracle language; the reverse computation rejects
otherwise. We call this a {\it reverse oracle call}.  
Although a reverse oracle call does not have the format of an oracle 
call, it can easily be simulated by an oracle call and a few more 
transitions. 
Hence we have:

\begin{pro} \label{Prop_invfPNP}
 \ Let $M$ be a polynomially balanced polynomial-time injective Turing
machine with {\sf NP}-oracle.
Then $M$ computes a polynomially balanced injective function.

The machine $M'$, obtained by running $M$ in reverse, is equivalent to a
polynomial-time injective Turing machines with {\sf NP}-oracle, and computes
$f^{-1}$ (where $f$ is the injective function computed by $M$). 
 \ \ \  \ \ \ $\Box$
\end{pro}
Conversely, every injective oracle Turing machine $M$ can be simulated by a
reverse injective Turing machine with the same oracle. Indeed, let 
$M'$ be the reverse of $M$, and let $M'_1$ be an ordinary injective oracle
machine that simulates $M'$. Finally, let $(M'_1)'$ be the reverse of 
$M'_1$. Then $(M'_1)'$ is a reverse injective Turing machine, with the same
oracle as $M$, simulating $M$.

So, there is no intrinsic difference between injective oracle Turing 
machines and reverse injective oracle Turing machines.
We could generalize injective Turing machines with oracle so as to allow 
oracle calls and reverse oracle calls in the same machine; but by the above 
discussion, it doesn't matter much whether we use oracle calls, reverse 
oracle calls, or both.

\begin{cor} \label{invfPoracle_inv} 
 \ \ ${\sf invfP}^{\sf (NP)}$ is an inverse monoid.
\end{cor}
{\bf Proof.} The composite of two polynomially balanced injective functions, 
computed by polynomial-time injective Turing machine with {\sf NP}-oracle, 
is also polynomially balanced and injective. It is computed by a 
polynomial-time injective Turing machine with {\sf NP}-oracle, obtained by
chaining the two machines. The disjoint union of two sets in {\sf NP} is in
{\sf NP}; thus the combined machine can use a single {\sf NP}-oracle 
language.  Hence, ${\sf invfP}^{\sf (NP)}$ is a monoid. 
By Prop.\ \ref{Prop_invfPNP}, ${\sf invfP}^{\sf (NP)}$ is an inverse 
monoid.
 \ \ \ $\Box$

\bigskip

\noindent {\bf Remark.} Not all functions in ${\sf fP}^{\sf (NP)}$ have an 
inverse in ${\sf fP}^{\sf (NP)}$ (unless the polynomial hierarchy {\sf PH} 
collapses, see Sect.\ 6 of \cite{s1f}). Here, ${\sf fP}^{\sf (NP)}$ denotes
the set of polynomially balanced functions computed by deterministic
polynomial-time Turing machines with {\sf NP}-oracle.    
But all functions in ${\sf invfP}^{\sf (NP)}$ have an inverse in 
${\sf invfP}^{\sf (NP)}$ (by Cor.\ \ref{invfPoracle_inv}).

\begin{pro} \label{injVSdetandrev} 
 \ Suppose $f: A^* \to A^*$ is an injective and polynomially balanced 
function such that both $f$ and $f^{-1}$ can be computed by deterministic 
polynomial-time Turing machines with {\sf NP}-oracle. 
Then $f$ can be computed by an {\em injective} polynomial-time Turing 
machine with {\sf NP}-oracle.
\end{pro}
{\bf Proof.} Bennett's proof in \cite{Bennett73, Bennett89} applies to 
Turing machines with {\sf NP}-oracle. An injective Turing machine for $f$ 
(with {\sf NP}-oracle) is obtained by first using the Turing machine for $f$ 
(on input $x$), but with a history tape added (which makes the machine 
injective). The input $x$ is still present. Once $f(x)$ has been computed, 
a copy of it is made on the output tape. The history tape is then used to 
run the previous computation in reverse, thus erasing the history tape
and the work-tape copy of $f(x)$. The input $x$ is still on the input tape
at this moment, and a copy of $f(x)$ is on the output tape. See Lemma 1 of
\cite{Bennett89} (except that now the machine also has an {\sf NP}-oracle).

To erase $x$ injectively, the Turing machine for $f^{-1}$ is used on input 
$f(x)$ (copied from the output tape). A history tape is added to make the
computation injective; $x$ is (re-)computed (while $f(x)$ is kept on the
output tape). Once $x$ has been (re-)computed, the history tape is used to 
run the previous computation of the $f^{-1}$-machine in reverse, thus erasing 
the history tape, as well as $x$. 
A copy of $f(x)$ is kept as the output tape.
See Theorem 2(b) of \cite{Bennett89} (except that now the machine also has 
an {\sf NP}-oracle).
 \ \ \ $\Box$

\medskip

We will prove next that $\, {\sf invfP}^{\sf (NP)}$ is finitely generated.
For this we want to replace all {\sf NP}-oracles by one {\sf NP}-complete 
set. Moreover, we want the reduction functions to be in {\sf invfP}.

Recall that for sets $L_1, L_2 \subseteq A^*$, a {\em many-one reduction}
from $L_1$ to $L_2$ is a polynomial-time computable function $f$: 
$A^* \to A^*$ such that $L_1 = f^{-1}(L_2)$;  equivalently, 
$(\forall x \in A^*)[ \, x \in L_1 \Leftrightarrow f(x) \in L_2 \, ]$.
A {\em one-one reduction} is a many-one reduction $f$ that is injective.

\begin{defn} \label{invfPRed}
 \ For sets $L_1, L_2 \subseteq A^*$, an {\em {\sf invfP}-reduction} from
$L_1$ to $L_2$ is any element $f \in {\sf invfP}$ such that 
$L_1 = f^{-1}(L_2)$.
\end{defn}
This definition generalizes the usual many-one reductions in the sense 
that the elements of {\sf invfP} are {\em partial} functions (with 
domain in {\sf P}). Note that if $L_1 = f^{-1}(L_2)$ then 
$L_1 \subseteq {\sf Dom}(f)$.

By the definition of {\sf invfP}, if $f \in {\sf invfP}$ then 
$f^{-1} \in {\sf invfP}$, and ${\sf Im}(f) \in {\sf P}$  (by Prop.\ 1.9 in
\cite{s1f}, since $f$ is regular).
If $f$ reduces $L_1$ to $L_2$ then $f^{-1}$ reduces $L_2 \cap {\sf Im}(f)$ 
to $L_1$. So, unless $L_2 \subseteq {\sf Im}(f)$, $f^{-1}$ does {\sl not} 
reduce $L_2$ to $L_1$, and {\sf invfP}-reducibility is not a symmetric 
relation.  

As we saw in Prop.\ \ref{invfPinjTM}, an {\sf invfP}-reduction can be 
computed by a polynomially balanced polynomial-time injective Turing 
machine.

\begin{pro} \label{NPcomplINVred}
 \ There exists languages that are {\sf NP}-complete with respect to 
{\sf invfP}-reductions.
\end{pro}
{\bf Proof.} An example is the ``universal {\sf NP}-complete language'' of
Hartmanis \cite{Hart}, defined as follows (slightly reformulated): 

\medskip

$L_{\rm univ}^{\sf NP} \ = \ $
$\{ {\sf code}(w) \ 11 \, x \, 11 \, 0^{|w| \cdot p_w(|x|)} :$

\smallskip

\hspace{0.7in} $w$ is a nondeterministic polynomial-time Turing machine program,

\hspace{0.7in}  $p_w$ is the built-in polynomial of $w$, and $x \in L_w \}$.

\medskip

\noindent In \cite{Hart} only many-one reductions were considered, but 
{\sf invfP}-reductions can be used too.
Indeed, for any language $L_v \in {\sf NP}$ accepted by a nondeterministic
polynomial-time Turing machine with program $v$, we define the function 
$g_v$ for all $x \in A^*$ by

\medskip

$x \ \longmapsto$
$ \ {\sf code}(v) \ 11 \, x \, 11 \, 0^{|v| \cdot p_v(|x|)}$. 

\medskip

\noindent Then $g_v$ belongs to {\sf invfP}, since a fixed program $v$ is 
chosen for the language $L_v$. And $g_v$ is a {\sf invfP}-reduction reduction 
from $L_v$ to $L_{\rm univ}^{\sf NP}$.  So,  $L_{\rm univ}^{\sf NP}$ is 
complete with respect to {\sf invfP}-reductions.
 \ \ \ $\Box$

\begin{pro} \label{invfPNP_FinGen}
 \ The inverse monoid $\, {\sf invfP}^{\sf (NP)}$ is finitely generated, 
as a monoid.
\end{pro}
{\bf Proof.} By Prop.\ \ref{NPcomplINVred} there exist sets that are 
{\sf NP}-complete with respect to {\sf invfP}-reductions.  Hence, for all 
polynomially balanced polynomial-time injective Turing machines with 
{\sf NP}-oracle, we can use a fixed set $N$ as the oracle, where $N$ is 
{\sf NP}-complete with respect to {\sf invfP}-reductions.
This changes the time-complexity of each function in 
${\sf invfP}^{\sf (NP)}$ by a polynomial amount.  

Now the proof of Prop.\ \ref{invfP_FinGen} goes through if we replace 
{\sf invfP}-programs by ${\sf invfP}^{\sf (NP)}$-programs.
The latter programs use the three additional states $q_{\sf qu}$, 
$q_{\sf yes}$, and $q_{\sf no}$, and one instruction for implementing a 
call to the oracle $N$.
We replace ${\sf injEv}_q$ by ${\sf injEv}_q^{(N)}$, where the latter is
computed by an injective Turing machine that is similar to the one for 
${\sf injEv}_q$, but with oracle calls to $N$ added.
 \ \ \ $\Box$

\bigskip

In the next section we will consider inverses and co-inverses that are 
computed by polynomially balanced polynomial-time injective Turing machines 
with {\sf NP}-oracle.

%%%%%%%%%%%%%%%%%%%%%%%%%%%%%%%%%%%%%%%%%%%%%%%%%%%%%%%%%%%%
\section{Inverses and co-inverses of any element of {\sf fP}}

If {\sf P} $\neq$ {\sf NP} then there exist non-regular functions 
$f \in {\sf fP}$. In that case we can nevertheless consider the injective
inverses and co-inverses of $f$ in ${\sf invfP}^{\sf (NP)}$, as we shall see
in Prop.\ \ref{fmin_in_invfPNP}.
First we derive some general results about injective inverses and 
co-inverses.

%%%%%%%%%%%%%
\subsection{Inverses, co-inverses, sub-inverses}

Here we consider arbitrary functions $f, f': A^* \to A^*$, unless more 
precise conditions are stated.
A subfunction of a function $f$ is, by definition, any function $g$ such that
$g \subseteq f$; equivalently, ${\sf Dom}(g) \subseteq {\sf Dom}(f)$ and for
all $x \in {\sf Dom}(g)$:  $g(x) = f(x)$.

\begin{lem} \label{InvDomIm}
 \ If $\, f f' f = f \,$ then $\, {\sf Im}(f) \subseteq {\sf Dom}(f')$.
\end{lem}
{\bf Proof.} If $y = f(x) \in {\sf Im}(f)$ then
$f f'(y) = f f' f(x) = f(x) = y$, so $f f'(y)$ is defined, hence $f'(y)$ is
defined, so $y \in {\sf Dom}(f')$.
  \ \ \ $\Box$

\begin{lem} \label{inj_co_invDomIm} \hspace{-0,13in} {\bf .}

\noindent {\bf (1)} \ If $\, f' f f' = f' \, $ then 
$\, {\sf Im}(f') \subseteq {\sf Dom}(f)$.

\noindent {\bf (2)} \ If in addition $f'$ is {\em injective}, then 
$\, {\sf Dom}(f') \subseteq {\sf Im}(f)$.
This does not hold in general when $f'$ is not injective.
\end{lem}
{\bf Proof.} (1) This is equivalent to Lemma \ref{InvDomIm}, up to notation. 

\noindent (2)  If $z \in {\sf Dom}(f')$ then $f'(z) = f' f f'(z)$ is 
defined, hence $f f'(z)$ is defined and $f f'(z) \in {\sf Im}(f)$.
When $f'$ is injective we can consider its set-theoretic inverse $f'^{-1}$, 
and we have $f'^{-1} f' = {\sf id}_{{\sf Dom}(f')}$.
Applying $f'^{-1}$ on the left to $f'(z) = f' f f'(z)$ yields
${\sf id}_{{\sf Dom}(f')} \circ f' f(z) = z$ (if $z \in {\sf Dom}(f')$).
Hence ${\sf id}_{{\sf Dom}(f')} \circ f f'(z)$ is defined, and 
${\sf id}_{{\sf Dom}(f')} \circ f f'(z) = f f'(z) = z$, 
so $z \in {\sf Im}(f)$.

When $f'$ is not injective, the result might not hold. Consider for example 
$f = \{(a,b)\}$ and $f' = \{(a, a), (b,a)\}$. Then $f'ff' =f' \,$ (and 
$f f' f = f$ as well, so $f$ and $f'$ are mutual inverses); but 
${\sf Dom}(f') = \{a, b\} \, \not\subseteq \, \{b\} = {\sf Im}(f)$.
 \ \ \ $\Box$

\begin{cor} \label{inj_mutinvinvDomIm}
 \ If $f, f'$ are mutual inverses and $f'$ is injective then 
${\sf Dom}(f') = {\sf Im}(f)$. 
This does not hold in general when $f'$ is not injective.
 \ \ \   \ \ \ $\Box$
\end{cor}

\begin{lem} \label{fvsfprimeinv}
 \ If $\, f' f f' = f' \, $ and $f'$ is {\em injective}, then
 \ $f'^{-1} \ = \ f|_{{\sf Im}(f')}$; \ equivalently,

\smallskip

 \ \ \  \ \ \         
$f'  \ = \ \big(f|_{{\sf Im}(f')}\big)^{-1}$
$ \ = \ f^{-1} \ \cap \ {\sf Dom}(f') \times {\sf Im}(f')$.

\medskip

\noindent Hence, $\, f'^{-1} \, \subseteq \, f$.
\end{lem}
{\bf Proof.} Applying $f'^{-1}$ on the left and the right to 
$f' = f' f f'$ yields 
$\, {\sf id}_{{\sf Dom}(f')} \circ f \circ {\sf id}_{{\sf Im}(f')}$ 
$ = f'^{-1}$.  This is equivalent to 
$\, f|_{{\sf Im}(f') \, \cap \, f^{-1}({\sf Dom}(f'))} = f'^{-1}$.
Hence, $f'^{-1} \subseteq f$, so $\, f' = (f'^{-1})^{-1} \subseteq f^{-1}$, 
where $f^{-1} = \{(y,x) : y = f(x)\}$.
Thus, ${\sf Im}(f') = f'({\sf Dom}(f')) \subseteq f^{-1}({\sf Dom}(f'))$,
hence $\, {\sf Im}(f') \, \cap \, f^{-1}({\sf Dom}(f')) = {\sf Im}(f')$; so
$\, f|_{{\sf Im}(f') \, \cap \, f^{-1}({\sf Dom}(f'))} = $
$f|_{{\sf Im}(f')} = f'^{-1}$.

The equality 
$\, f' \, = \, f^{-1} \, \cap \, {\sf Dom}(f') \times {\sf Im}(f')$
is just a logical reformulation of $f' = (f|_{{\sf Im}(f')})^{-1}$.
 \ \ \ $\Box$

\begin{cor} \label{COinvSubfunc}
 \ If $f'$ is an injective co-inverse of $f$, then $f'$ is a mutual inverse 
of $\, f|_{{\sf Im}(f')} \, $ (which is an injective subfunction of $f$).
Moreover, $f'$ is uniquely determined by $f$ and ${\sf Im}(f')$, as
$\, f' = (f|_{{\sf Im}(f')})^{-1}$. 
\end{cor}
{\bf Proof.} Since $\, f' = {\sf id}_{{\sf Im}(f')} \circ f'$, we have 
$\, f' = f' f f' = f' \circ f \circ {\sf id}_{{\sf Im}(f')} \circ f'$ 
$= f' \circ f|_{{\sf Im}(f')} \circ f'$.
Moreover, since $f|_{{\sf Im}(f')} = f'^{-1}$ (by Lemma \ref{fvsfprimeinv}),
we have $\, f|_{{\sf Im}(f')} = f'^{-1} = f'^{-1} f' f'^{-1} = $
$f|_{{\sf Im}(f')} \circ f' \circ f|_{{\sf Im}(f')}$.
 \ \ \ $\Box$

\begin{cor} \label{COinvTOinv}
 \ If $\, f' f f' = f' \,$ and $\, {\sf Dom}(f') = {\sf Im}(f)$, then 
$f'$ is an injective mutual inverse of $f$.
\end{cor}
{\bf Proof.} When ${\sf Dom}(f') = {\sf Im}(f)$ then ${\sf Im}(f')$ is a
choice set of $f$. Indeed, for every $y \in {\sf Im}(f)$, $f'(y)$ is defined
and $f'(y) \in f^{-1}(y) \,$ (since $f' \subset f^{-1}$ by the second 
equality in Lemma \ref{fvsfprimeinv}).  
Hence $f'$ is an inverse of $f$.
Injectiveness follows from the observations made in the Introduction and at 
the beginning of Section 2.  
 \ \ \ $\Box$

\begin{thm} \label{COinvOFprod} 
{\bf (anti-homomorphic property of co-inverses).}

\noindent Suppose $f_i'$ is an injective co-inverse of $f_i$ for $i = 1, 2$. 
Then $\, f_1' f_2' \,$ is an injective co-inverse of $\, f_2 f_1$. 
\end{thm}
{\bf Proof.} We first observe that for subidentities ${\sf id}_X$,
${\sf id}_Y$, we have $\,  {\sf id}_X \ {\sf id}_Y = {\sf id}_{X \cap Y}$
$= {\sf id}_Y \ {\sf id}_X$. 
Since $f_i'$ is injective, 
$f_i' f_i'^{-1} = {\sf id}_{{\sf Im}(f_i')}$, and 
$\, f_i'^{-1} f_i' = {\sf id}_{{\sf Dom}(f_i')}$.  Hence 

\smallskip

$\, f_1' f_2' f_2'^{-1} f_1'^{-1} f_1' f_2' \, = \, $
$f_1' \ {\sf id}_{{\sf Im}(f_2')} \ {\sf id}_{{\sf Dom}(f_1')} \ f_2' \,$
$ = \, f_1' \ {\sf id}_{{\sf Dom}(f_1')} \ {\sf id}_{{\sf Im}(f_2')} \ f_2'$
$ \, = \, f_1' f_2'$. 

\smallskip

\noindent So, 
$f_1' f_2' \, = \, f_1' f_2' f_2'^{-1} f_1'^{-1} f_1' f_2' \ \subseteq \ $
$f_1' f_2' f_2 f_1  f_1' f_2' \,$ (by Lemma \ref{fvsfprimeinv}). 
Now, since $f_1' f_2' \, \subseteq \, f_1' f_2' f_2 f_1  f_1' f_2'$, we have
$\, f_1' f_2'(x) = f_1' f_2' f_2 f_1  f_1' f_2'(x) \, $ when 
$f_1' f_2'(x)$ is defined.  And when $f_1' f_2'(x)$ is not defined, 
$f_1' f_2' f_2 f_1  f_1' f_2'(x)$ is not defined either.
Thus, $f_1' f_2' \, = \, f_1' f_2' f_2 f_1  f_1' f_2'$.  
 \ \ \ $\Box$

\bigskip

\noindent {\bf Remark.} Theorem \ref{COinvOFprod}, though straightforward,
is remarkable as it applies to injective {\it co}-inverses but not to injective 
inverses or mutual inverses. 

The latter can be illustrated by the following example:
Let $f_1 = \{(0,0), (1,0)\} = f_2$, and let $f_1' = \{(0,1)\} = f_2'$, 
so $f_i$ and $f_i'$ are mutual inverses $\, (i = 1, 2)$, and $f_i'$ is
injective.
But $f_1' \circ f_2' = \theta \,$ (the empty map), which is not an inverse
of $f_2 f_1 = \{(0,0), (1,0)\}$. 

\bigskip

In \cite{s1f} (Prop.\ 6.1) we showed that every element $f \in {\sf fP}$ 
has an inverse in ${\sf fP}^{\sf (NP)}$; one such inverse is $f_{\sf min}'$, 
defined by 

\smallskip

 \ \ \ \ \ \ $f_{\sf min}'(y) \ = \ \left\{
   \begin{array}{ll}
    {\sf min}(f^{-1}(y))    & \ \ \ \mbox{if $y \in {\sf Im}(f)$, } \\
    {\rm undefined}         & \ \ \ \mbox{otherwise,}
   \end{array}  \right.
$

\smallskip

\noindent where the {\sf min} operation is taken with respect to the 
length-lexicographic order, $\le_{\sf llex}$, of $\{0,1\}^*$. Recall that by 
definition, $u \le_{\sf llex} v \,$ iff $\, |u| < |v|$ or $[ \, |u| = |v|$ 
and $u \le_{\sf dict} v \,]$; here, $\le_{\sf dict}$ is the dictionary 
order of $\{0,1\}^*$ determined by $0 <_{\sf dict} 1$.

\begin{pro} \label{fmin_properties}
For every $f \in {\sf fP}$ we have: \  

\noindent {\bf (1)} \ $f \circ f_{\sf min}' = {\sf id}_{{\sf Im}(f)}$;
 \ hence, $f \circ f_{\sf min}' \circ f = f$; 
 \ also \ $f_{\sf min}' \circ f \circ f_{\sf min}' = f_{\sf min}'$;
 \ \ {\bf (2)} \ $f_{\sf min}'$ is injective;  

\smallskip

\noindent 
{\bf (3)} \ ${\sf Im}(f_{\sf min}') \in {\sf coNP}$;
 \ \ {\bf (4)} \ ${\sf Dom}(f_{\sf min}') = {\sf Im}(f) \in {\sf NP}$;
 \ \ {\bf (5)} \ $f_{\sf min}'^{-1} \in {\sf fP}^{\sf (NP)}$. 
\end{pro}
{\bf Proof.} \ (1) is obvious from the definition of $f_{\sf min}'$.

\noindent (2) Injectiveness of $f_{\sf min}'$ follows from the fact 
that the sets $f^{-1}(y)$ (for $y \in {\sf Im}(f)$) are two-by-two disjoint.
 
\noindent (3) We have:
$\, x \in {\sf Im}(f_{\sf min}')$ \ iff \ $x \in {\sf Dom}(f)$ and 
$(\forall z)[ \, z <_{\sf llex} x$ $\Rightarrow f(z) \neq f(x) \, ] \,$ 
(the latter sentence expresses that $x$ is the $\le_{\sf llex}$-minimum 
element in its ${\sf mod}_f$-class).  
Since ${\sf Dom}(f) \in {\sf P}$, since the relation
$\, \{ (z,x) : z <_{\sf llex} x \Rightarrow f(z) \neq f(x)\} \,$ is in 
{\sf P}, and since this relation is universally quantified (by
$(\forall z, z <_{\sf llex} x)$),  it follows that 
${\sf Im}(f_{\sf min}') \in {\sf coNP}$. 

\noindent (4) By the definition of $f_{\sf min}'$ we have 
${\sf Dom}(f_{\sf min}') = {\sf Im}(f)$. And ${\sf Im}(f) \in {\sf NP}$ for
every $f \in {\sf fP}$. 

\noindent (5) As a consequence of (2), 
$f_{\sf min}'^{-1} \in {\sf fP}^{\sf (NP)}$. Indeed,
$f_{\sf min}'^{-1} = f|_{{\sf Im}(f_{\sf min}')}$ (by Lemma 
\ref{fvsfprimeinv}), and $f \in {\sf fP}$, and 
${\sf Im}(f_{\sf min}') \in {\sf coNP}$, so $f_{\sf min}'^{-1}$ is computed 
by any Turing machine for $f$ with an added oracle call to 
${\sf Im}(f_{\sf min}')$. 
 \ \ \ $\Box$

\medskip
 
In the proof of Prop.\ 6.1 in \cite{s1f} we showed that $f_{\sf min}' \in$
${\sf fP}^{\sf (NP)}$. Next we strengthen this by showing that the oracle 
Turing machine that computes $f_{\sf min}'$ is an {\em injective} Turing 
machine with an {\sf NP}-oracle.  
In Def.\ \ref{DEFinvfPNP} we introduced the set ${\sf invfP}^{\sf (NP)}$ of 
polynomially balanced functions computed by polynomial-time injective 
Turing machines with an {\sf NP}-oracle.

\begin{pro} \label{fmin_in_invfPNP}
 \ Every $f \in {\sf fP}$ has an inverse $f'$ in 
$\, {\sf invfP}^{\sf (NP)} \,$ (hence $f'^{-1} \in {\sf invfP}^{\sf (NP)}$), 
with the additional property that $\, {\sf Im}(f') \in {\sf coNP}$.
\end{pro}
{\bf Proof.} \ Let $f'$ be $f_{\sf min}'$, which is injective as we just 
saw, hence $f'^{-1}$ is an injective function.  
Both $f'$ and $f'^{-1}$ are in ${\sf fP}^{\sf (NP)}$; this holds for $f'$ 
by Prop.\ 6.1 in \cite{s1f}, and for $f'^{-1}$ by Prop.\ 
\ref{fmin_properties}. 
Hence by Prop.\ \ref{injVSdetandrev},  
$f' \in {\sf invfP}^{\sf (NP)}$.

Since $f'$ is computed by an injective Turing machine (with oracle), 
$f'^{-1}$ can be computed by the same machine run backwards; so
$f'^{-1} \in {\sf invfP}^{\sf (NP)}$.

We proved in Prop.\ \ref{fmin_in_invfPNP} that 
${\sf Im}(f_{\sf min}') \in {\sf coNP}$.
 \ \ \ $\Box$

\medskip

\noindent {\bf Remark.} As an alternative proof of Prop.\ 
\ref{fmin_in_invfPNP}, observe that the algorithm for computing 
$f_{\sf min}'(y)$ in the proof of Prop.\ 6.1 in \cite{s1f} actually 
describes an {\em injective} Turing machine; it makes oracle calls to the 
set $\, \{(z,u) : z \in f(u A^*)\} \,$ (which is in {\sf NP} for any 
fixed $f \in {\sf fP}$).

\bigskip

\noindent Based on the anti-homomorphic property of co-inverses (Theorem
\ref{COinvOFprod}) we define:

\begin{defn} \label{DEFcofP}
The {\em monoid of co-inverses} of {\sf fP} is

\medskip

 \ \ \  \ \ \ ${\sf cofP} \ =$ 
$ \ \{ f' \in {\sf invfP}^{\sf (NP)} : \, $
 $f'$ is a co-inverse of some element of ${\sf fP} \}$.  

\medskip

\noindent Equivalently, an element $f'$ of ${\sf invfP}^{\sf (NP)}$ belongs 
to ${\sf cofP}$ \ iff \ $f'$ has an inverse in {\sf fP}.
\end{defn}
By Theorem \ref{COinvOFprod}, ${\sf cofP}$ is indeed a {\em monoid}.  

Also, $f_{\sf min}' \in {\sf cofP}$ for every $f \in {\sf fP}$, 
since $f$ is a mutual inverse of $f_{\sf min}'$ (Prop.\ 
\ref{fmin_properties}). Hence, for every $f \in {\sf fP}$, 
$\, {\sf cofP}$ contains a mutual inverse of $f$.

\medskip

One could also consider the {\em monoid of inverses} of {\sf fP}, namely
$\, \{g \in {\sf invfP}^{\sf (NP)} : \, g$ is an inverse of some
element of ${\sf fP} \}$.
However, this is just all of ${\sf invfP}^{\sf (NP)}$, since every element of
${\sf invfP}^{\sf (NP)}$ is an inverse of the empty map $\theta \in {\sf fP}$.
So, {\em being an inverse} of an element in {\sf fP} is a trivial notion;
but {\em having an inverse} in {\sf fP} (i.e., being a co-inverse of an
element in {\sf fP}) is non-trivial.
The question whether ${\sf cofP}$ is equal to
${\sf invfP}^{\sf (NP)}$ has an interesting ``answer'', which is presented
in Prop.\ \ref{cofPvsfPNP} below.

\begin{pro} \label{cofP_cap_fP}
 \ \ \ ${\sf cofP} \ \cap \ {\sf fP} \ = \ {\sf invfP}$. 
\end{pro}
{\bf Proof.} $[\supseteq]$ \ If $f' \in {\sf invfP}$ then
$f'^{-1} \in {\sf invfP} \subseteq {\sf fP}$. Moreover,
$f' \, f'^{-1} \, f' = f'$, so $f'$ is a co-inverse of an element in
{\sf fP}.
Hence $f' \in {\sf cofP}$.

\noindent $[\subseteq]$ \ If $f' \in {\sf cofP} \cap {\sf fP}$
then there exists $f \in {\sf fP}$ such that $f' f f' = f'$, hence $f'$ is
a regular element of {\sf fP}. Since we assume
$f' \in {\sf cofP}$, it follows that $f'$ is injective. So, $f'$
is an injective regular element of {\sf fP}, i.e., $f' \in {\sf invfP}$.
 \ \ \ $\Box$

\bigskip

On the other hand, ${\sf invfP}^{\sf (NP)} \ \cap \ {\sf fP}$
contains all injective one-way functions. In fact,
${\sf invfP}^{\sf (NP)} \, \cap \, {\sf fP} \ = \ {\sf invfP}$ \ iff
 \ injective one-way functions do not exist. The latter is the case iff
$\, {\sf P} = {\sf UP} \,$ (see e.g.\ section 12.1 in \cite{Papadim}).
But {\sf cofP} contains no one-way function, since every element of 
{\sf cofP} has an inverse in {\sf fP}.

\begin{defn} \label{DEFinjfP}
 \ Let {\sf injfP} denote the monoid of all injective functions in {\sf fP}.
\end{defn}
The set of regular elements of {\sf injfP} is exactly the submonoid 
{\sf invfP} (by the definition of {\sf invfP}), hence the set of injective 
one-way functions (for worst-case complexity, see Prop.\ \ref{existsInj1wf}) 
is $\, {\sf injfP} - {\sf invfP}$.
We saw that {\sf invfP} and ${\sf invfP}^{\sf (NP)}$ are regular monoids (in
fact, inverse monoids).
On the other hand, {\sf injfP} is regular iff injective one-way functions do 
not exist. Thus we have:
 
\smallskip

 {\it injective one-way functions exist \ iff } 
 \ \ ${\sf invfP} \neq {\sf injfP} \neq {\sf invfP}^{\sf (NP)}$.

\smallskip

\noindent Since {\sf invfP} consists of the regular elements of {\sf injfP}, 
we also have:

\smallskip

 {\it injective one-way functions exist \ iff } 
 \ \ ${\sf invfP} \neq {\sf injfP}$.

\smallskip

\noindent Since ${\sf invfP}^{\sf (NP)}$ is regular we have,  

\smallskip

 {\it if injective one-way functions exist, then }
 \ ${\sf injfP} \neq {\sf invfP}^{\sf (NP)}$. 

\smallskip

\noindent We will also see (Theorem \ref{cofP_vs_P_NP}) that one-way 
functions exist \ iff \ ${\sf invfP} \neq {\sf invfP}^{\sf (NP)}$.
 
We saw that {\sf cofP} contains no one-way functions (see the
observations after Prop.\ \ref{cofP_cap_fP}), so we have:
  \ $\, {\sf injfP} \ \cap \ {\sf cofP} \ = \ {\sf invfP}$.
 \ Moreover:

\begin{pro} \label{invfPNPcapfP}
 \ ${\sf invfP}^{\sf (NP)} \, \cap \, {\sf fP} \ = \ {\sf injfP}$.
\end{pro}
{\bf Proof.} $[\subseteq]$ \ Every $f \in {\sf invfP}^{\sf (NP)}$ is
injective, so the inclusion $\subseteq$ holds. 

\noindent $[\supseteq]$ \ If $f \in {\sf injfP}$ then 
$f_{\sf min}' = f^{-1} \in {\sf invfP}^{\sf (NP)}$, by Prop.\
\ref{fmin_in_invfPNP}. It follows that $(f^{-1})^{-1} = f \in$
${\sf invfP}^{\sf (NP)}$, by  Prop.\ \ref{Prop_invfPNP}.  
Hence $f \in {\sf invfP}^{\sf (NP)} \cap {\sf injfP}$ $\subseteq$
${\sf invfP}^{\sf (NP)} \cap {\sf fP}$.  
 \ \ \ $\Box$

\bigskip

\bigskip

\bigskip

\bigskip

%%%%%%%%%% Fig. 1 %%%%%%%%%%%%%%%%%%%%%%%%%%%%%%%%%%%%%%%%%%%%%%%%
\unitlength=.8mm 
\begin{picture}(80,110)

%% invfP^(NP)
\put(45,60){\oval(85,75)}
\put(18,88){\makebox(0,0)[cc]{${\sf invfP}^{\sf (NP)}$}}

%% cofP
\put(55,50){\oval(50,40)}
\put(42,63){\makebox(0,0)[cc]{\sf cofP}}

%% invfP
\put(70,50){\oval(19,39)[r]}
\put(70,38){\makebox(0,0)[cc]{\sf invfP}}
\put(60.5, 31){\line(0,1){38}}
\put(60.5, 30.5){\line(1,0){10}}
\put(60.5, 69.5){\line(1,0){10}}

%% injfP
\put(73,83){\makebox(0,0)[cc]{\sf injfP}}

%% fP
\put(90,64){\oval(60,97)}
\put(110,101){\makebox(0,0)[cc]{\sf fP}}

\end{picture}

\vspace{-5mm}

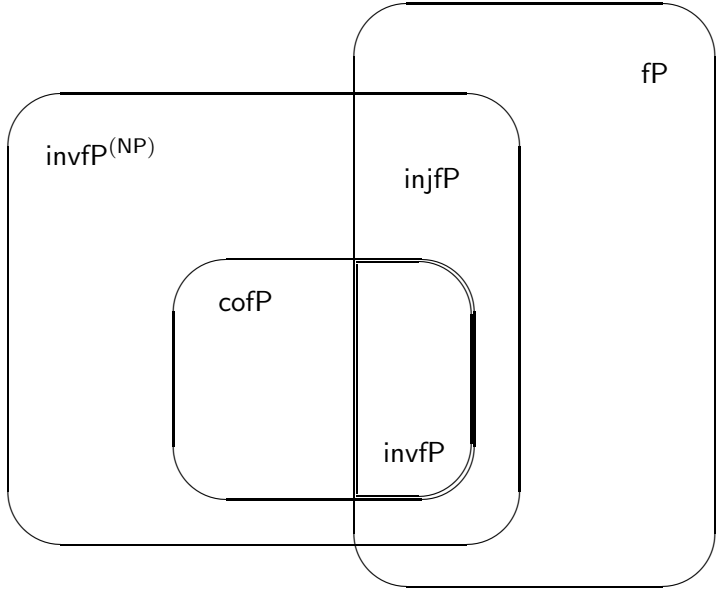
{\small \sl Figure 1: \ Relations between $\, {\sf fP}$, 
$\, {\sf invfP}^{\sf (NP)}$, $\, {\sf cofP}$, $\, {\sf invfP}$, and 
$\, {\sf injfP}$. }

\bigskip

\bigskip

\bigskip

\noindent There is an interesting connection between {\em co-inverses} 
and {\em sub-inverses} of $f$.

\begin{defn} \label{Subinv_def}
 \ A function $g'$ is a {\em sub-inverse} of a function $f$ iff there exists
$g \subseteq f$ such that $g g' g = g$ and $g' g g' = g'$.
I.e., the sub-inverses of $f$ are the mutual inverses of the subfunctions
of $f$.
\end{defn}

\begin{pro} \label{subVSco}
 \ For any function $f: A^* \to A^*$, the injective co-inverses of $f$ are
the same as the injective sub-inverses of $f$. 
\end{pro} 
{\bf Proof.} $[\Rightarrow]$ Let $f'$ be an injective co-inverse of $f$;
then ${\sf Dom}(f') \subseteq {\sf Im}(f)$ and 
${\sf Im}(f') \subseteq {\sf Dom}(f)$, by Lemmas \ref{InvDomIm} and
\ref{inj_co_invDomIm}. 
Let us consider the following restriction of $f$: \ $g = f|_{{\sf Im}(f')}$.
We claim that $f'$ is a mutual inverse of $g$, hence $f'$ is a sub-inverse 
of $f$.
Recall that $f' = (f|_{{\sf Im}(f')})^{-1}$, by Lemma \ref{fvsfprimeinv}.
It follows from this that 
$\, g f' g = f|_{{\sf Im}(f')} \circ (f|_{{\sf Im}(f')})^{-1}$ $\circ$ 
$f|_{{\sf Im}(f')} = f|_{{\sf Im}(f')} = g$; similarly it follows that 
$\, f' g f' = (f|_{{\sf Im}(f')})^{-1}$ $\circ$ 
$f|_{{\sf Im}(f')}$ $\circ$ $(f|_{{\sf Im}(f')})^{-1}$
$ = (f|_{{\sf Im}(f')})^{-1} = f'$. 

So, $f'$ is an injective sub-inverses of $f$.

\smallskip

\noindent $[\Leftarrow]$ Let $g \subset f$ be any restriction of $f$, and 
let $g'$ be any injective mutual inverse of $g$; in other words, $g'$ is an
injective sub-inverse of $f$. 
We claim that $g'$ is a co-inverse of $f$. 

By Corollary \ref{inj_mutinvinvDomIm}, and Lemma \ref{fvsfprimeinv}, we
have: ${\sf Dom}(f') = {\sf Im}(f)$, and 
$g' = (g|_{{\sf Im}(g')})^{-1}$.
Since $g \subset f$ the latter implies $g' =  (f|_{{\sf Im}(f')})^{-1}$.
Hence, $g' f g' = (f|_{{\sf Im}(f')})^{-1} \circ f \circ$
$(f|_{{\sf Im}(f')})^{-1} = (f|_{{\sf Im}(f')})^{-1} = g'$.
So, $g'$ is a co-inverse of $f$.
 \ \ \ $\Box$

\medskip

\noindent Note that in part $[\Leftarrow]$ of the above proof, $g$ need not
be injective, but $g|_{{\sf Im}(g')} = f|_{{\sf Im}(g')}$ is injective.
Thus we have: 

{\it If $g'$ is an injective sub-inverse of $f$ and it is a mutual inverse 
of a subfunction of $f$, then $g'$ is also a mutual inverse of an} 
injective {\it subfunction of $f$.} 

\bigskip

\noindent The following gives a relation between two notions of sub-inverse.
We defined a sub-inverse of $f$ to be a mutual inverse of a subfunction of
$f$ (Def.\ \ref{Subinv_def}). We could, instead, have defined a sub-inverse 
of $f$ to be a subfunction of a mutual inverse of $f$.

\begin{pro} \label{SubofInv} \hspace{-0.13in} {\bf .}

\noindent {\bf (1)} \ If $g'$ is a subfunction of a  mutual inverse of $f$, 
then $g'$ is a sub-inverse of $f$. (For this fact, $g'$ need not be 
injective, and $f$ need not be in {\sf fP}.)

\smallskip

\noindent {\bf (2)} \ If $g' \in {\sf invfP}^{\sf (NP)}$ is a sub-inverse of 
$f \in {\sf fP}$, then $g'$ is a subfunction of some mutual inverse 
$f' \in {\sf invfP}^{\sf (NP)}$ of $f$.  
\end{pro}
{\bf Proof.} (1) Let $g' \subseteq f'$ with $f f' f = f$ and $f' f f' = f'$. 
Let $g = f g' f$; then $g \subseteq f f' f = f$. 
And for all $x \in {\sf Dom}(g)$ we have 
$g g' g(x) = f f' f(x) = f(x) = g(x)$ (the latter since 
$x \in {\sf Dom}(g)$). So $g g' g = g$.
Moreover, for all $y \in {\sf Dom}(g')$ we have $g' g g'(y) = f' f f'(y)$,
since $g' \subseteq f'$ and $g \subseteq f$; and $f' f f'(y) = y$.
So, $g' g g' = g'$.

(2) Since $g' \in {\sf invfP}^{\sf (NP)}$, $g'$ is injective, hence it is an
injective co-inverse of $f$ (by Prop.\ \ref{subVSco}). Hence,
$g' = (f|_{{\sf Im}(g')})^{-1}$, and $g'$ is a mutual inverse of
$g = f|_{{\sf Im}(g')} \subseteq f$. 
We can extend $g'$ to the following mutual inverse of $f$:
For all $y \in {\sf Im}(f)$ let 

\smallskip

 \ \ \ \ \ \ $f'(y) \ = \ \left\{
   \begin{array}{ll}
    g'(y)                   & \ \ \ \mbox{if $\, y \in {\sf Dom}(g')$,} \\  
    f_{\sf min}'(y)         & \ \ \ \mbox{otherwise.}
   \end{array}  \right.
$

\smallskip

\noindent It follows that $f'$ is injective. Indeed, both $g'$ and 
$f_{\sf min}'$ are injective. Moreover, $g'(y_1) \ne f_{\sf min}'(y_2)$ for 
all $y_1 \in {\sf Dom}(g')$ and all $y_2 \in {\sf Im}(f) - {\sf Dom}(g')$,
since if we had equality then $f g'(y_1) = f f_{\sf min}'(y_2)$; but
$f g'(y_1) = y_1$ since $g' = (f|_{{\sf Im}(g')})^{-1}$, and  
$f f_{\sf min}'(y_2) = y_2$; so we would have $y_1 = y_2$, but 
$y_1 \in {\sf Dom}(g')$ and $y_2 \not\in {\sf Dom}(g')$.

Also, $f'$ is clearly a mutual inverse of $f$. 

Since $g'$ is computed by a injective Turing machine with 
{\sf NP}-oracle, and likewise $f_{\sf min}'$, there is an injective
Turing machine with {\sf NP}-oracle for $f'$; the Turing machine for $g'$,
being deterministic, can check whether $g'(y)$ is defined, i.e., whether
$y \in {\sf Dom}(g')$. Hence, $f' \in {\sf invfP}^{\sf (NP)}$.
 \ \ \ $\Box$

%%%%%%%%%%%
\subsection{Connections with {\sf NP} }

The next theorem motivates the study of the monoid ${\sf cofP}$ (and 
of {\sf invfP}) in the context of the {\sf P} versus {\sf NP} problem.  
We prove an easy Lemma first.

\begin{lem} \label{regExt}
 \ Let $f, h$ be functions $A^* \to A^*$ and let $Z \subset A^*$.
If $\  f \circ h|_Z \circ f = f$, then $f h f = f$.
\end{lem}
{\bf Proof.} For $x \not\in {\sf Dom}(f)$, $f h f(x)$ and $f(x)$ are both
undefined. 
For $x \in {\sf Dom}(f)$, $\, h|_Z(f(x))$ is defined, otherwise
$f \circ h|_Z \circ f(x)$ would be undefined, and thus would not be equal
to $f(x)$. Hence we have  $f(x) \in Z$, so ${\sf Im}(f) \subseteq Z$.
It follows that $h|_Z(f(x)) = h(f(x))$, for all $x \in {\sf Dom}(f)$.
Hence $h|_Z \circ f = h \circ f$, thus
$f \circ h|_Z \circ f = f \, h \, f = f$.
 \ \ \ $\Box$

\begin{thm} \label{cofP_vs_P_NP}
 \ The following are equivalent: \ \ 

 \ ${\sf P} = {\sf NP}$; 

 \ ${\sf cofP} = {\sf invfP}$;

 \ ${\sf cofP} \subset {\sf fP}$; 

 \ ${\sf cofP} \, $ is regular;

 \ ${\sf invfP}^{\sf (NP)} \ = \ {\sf invfP}$.  
\end{thm}
{\bf Proof.} We use the following logical fact: 
If $X \Rightarrow Y_i$ and $\overline{X} \Rightarrow \overline{Y}_i$, 
then $X \Leftrightarrow Y_i$ (for $i = 1, 2, \ldots \, $).
Here, ``${\sf P} = {\sf NP}$'' plays the role of $X$.

\smallskip

\noindent $[\Rightarrow]$ \ We assume ${\sf P} = {\sf NP}$.
We obviously have ${\sf invfP} \subseteq {\sf cofP}$.
If ${\sf P} = {\sf NP}$ then ${\sf fP}^{\sf (NP)} = {\sf fP}$,
hence ${\sf cofP} \subseteq {\sf fP}$. Also, {\sf fP} is 
regular when ${\sf P} = {\sf NP}$.
Since all elements of ${\sf cofP}$ are injective, all 
elements of ${\sf cofP}$ are injective regular elements of {\sf fP}, 
i.e., ${\sf cofP} \subseteq {\sf invfP}$. In conclusion, 
${\sf cofP} = {\sf invfP}$.

If ${\sf cofP} = {\sf invfP}$ then ${\sf cofP} \subset {\sf fP}$.
The inclusion is strict because {\sf fP} has some non-injective elements. 

And if ${\sf  cofP} = {\sf invfP}$ then, since ${\sf invfP}$ is 
regular, ${\sf cofP}$ is regular.

\smallskip

\noindent $[\Leftarrow]$ \ If ${\sf P} \neq {\sf NP}$ then there exists 
$f \in {\sf fP}$ such that $f$ is {\em not} regular. For $f$ there exists a 
mutual inverse $f' \in {\sf invfP}^{\sf (NP)}$ such that 
${\sf Dom}(f') = {\sf Im}(f)$; hence, $f' \in {\sf cofP}$. 
Moreover, ${\sf Dom}(f') = {\sf Im}(f) \in {\sf NP}$; and 
$f'$ can be chosen so that in addition, ${\sf Im}(f') \in {\sf coNP}$ 
(e.g., by choosing $f' = f'_{\sf min}$ as in Prop.\ \ref{fmin_in_invfPNP}).
Hence, $f'^{-1} \in {\sf invfP}^{\sf (NP)}$ (since $f'^{-1} \, = \, $
$f|_{{\sf Im}(f')}$ and ${\sf Im}(f') \in {\sf coNP}$).
However, $f' \not\in {\sf fP}$, since  $f$ is not regular in {\sf fP}.

Let us assume by contradiction that the monoid ${\sf cofP}$ is regular. 
Then $f'$ (as above) is a regular element of ${\sf cofP}$, hence (by Lemma 
\ref{fmin_in_invfPNP}), $f'^{-1} \in {\sf cofP}$.  This 
implies (by the definition of ${\sf cofP}$) that 
$\, f'^{-1} \, h \, f'^{-1} = f'^{-1}$ for some $h \in {\sf fP}$.
Composing on the left and on the right with $f'$ yields
 \ ${\sf id}_{{\sf Im}(f')} \ h \ {\sf id}_{{\sf Dom}(f')} = f'$;
equivalently, $h|_{{\sf Dom}(f') \, \cap \, h^{-1}({\sf Im}(f'))} = f'$. 
Hence, $f f' f = $
$f \ h|_{{\sf Dom}(f') \, \cap \, h^{-1}({\sf Im}(f'))} \ f$ $=$ $f$.
It follows (by Lemma \ref{regExt}) that $h$ is an inverse of $f$, since 
the restriction $h|_{{\sf Dom}(f') \, \cap \, h^{-1}({\sf Im}(f'))}$ is 
an inverse of $f$.  Hence (since $f$ is not regular),
$h \not\in {\sf fP}$; but this contradicts the previous assumption that 
$h \in {\sf fP}$. 
Thus ${\sf cofP}$ is not regular.

Since the monoid {\sf invfP} is regular, we conclude that
${\sf cofP} \neq {\sf invfP}$. 

And ${\sf cofP} \neq {\sf invfP}$ implies 
${\sf cofP} \not\subset {\sf fP}$; indeed, if we had 
${\sf cofP} \subset {\sf fP}$ then we would have 
${\sf cofP} = {\sf cofP} \, \cap \, {\sf fP} =$ 
${\sf invfP}$. 

Finally, if ${\sf invfP}^{\sf (NP)} = {\sf invfP}$ then 
${\sf cofP} = {\sf invfP}$, hence by what we proved above, 
${\sf P} = {\sf NP}$. Conversely, if  ${\sf P} = {\sf NP}$ then every 
polynomial-time injective Turing machine with {\sf NP}-oracle is equivalent
to a polynomial-time deterministic Turing machine (without oracle). Thus
if $f \in {\sf invfP}^{\sf (NP)}$ (hence by Prop.\ \ref{Prop_invfPNP}, 
$f^{-1} \in {\sf invfP}^{\sf (NP)}$), then $f, f^{-1} \in {\sf fP}$. 
By Bennett's results \cite{Bennett73, Bennett89} (especially Lemma 1 and 
Theorem 2 (b) in \cite{Bennett89}), this implies $f \in {\sf invfP}$.
  \ \ \ $\Box$

\bigskip

Although every element of ${\sf cofP}$ has an inverse in 
{\sf fP} (by the definition of ${\sf cofP}$), not every element 
of ${\sf cofP}$ has an inverse in ${\sf cofP}$, 
unless ${\sf P} = {\sf NP}$; the latter follows from the fact that 
${\sf cofP}$ is regular iff ${\sf P} = {\sf NP}$.

\bigskip

By definition, 
${\sf cofP} \subseteq {\sf invfP}^{\sf (NP)}$. We will now address the 
question whether this inclusion is strict; it turns out, interestingly, 
that this is equivalent to the question whether ${\sf P} \neq {\sf NP}$.

\begin{lem} \label{fprime_invINfPNP}
 \ Let $f \in {\sf fP}$, and let $f' \in {\sf cofP}$ be a co-inverse of 
$f$ such that $\, {\sf Im}(f') \in {\sf coNP}$.  Then
$\, f'^{-1} \in {\sf fP}^{\sf (NP)}$.

For every $f \in {\sf fP}$ there exists a co-inverse
$f' \in {\sf cofP}$ of $f$ such that
$\, {\sf Im}(f') \in {\sf coNP}$ and $\, {\sf Dom}(f') \in {\sf NP}$.
\end{lem}
{\bf Proof.} Since $f \in {\sf fP}$ and ${\sf Im}(f') \in {\sf coNP}$,
$f'^{-1} = f|_{{\sf Im}(f')}$ can be computed by a Turing machine for $f$
with {\sf NP}-oracle ${\sf Im}(f')$.

We saw above (and in \cite{s1f}) that for every $f \in {\sf fP}$ there
exists a mutual inverse (hence a co-inverse) $f'$ for which
${\sf Im}(f') \in {\sf coNP}$ and ${\sf Dom}(f') \in {\sf NP}$; e.g., let
$f' = f_{\sf min}'$.
 \ \ \ $\Box$

\begin{pro} \label{cofPvsfPNP}
 \ \ ${\sf P} \neq {\sf NP}$ \ iff 
 \ $\, {\sf cofP} \neq {\sf invfP}^{\sf (NP)}$.

\smallskip

In any case (i.e., whether ${\sf P} \ne {\sf NP}$ or not), ${\sf cofP}$ 
and ${\sf invfP}^{\sf (NP)}$ have the same set of idempotents.
\end{pro}
{\bf Proof.} $[\Leftarrow]$ \ If 
${\sf cofP} = {\sf invfP}^{\sf (NP)}$, then for every 
$f' \in {\sf invfP}^{\sf (NP)}$ there exists 
$f \in {\sf fP}$ such that $f' f f' = f'$. By Lemma \ref{fvsfprimeinv}, 
$f'^{-1} \subseteq f$, and by Lemma \ref{fprime_invINfPNP}, 
$f'^{-1} \in {\sf invfP}^{\sf (NP)}$. 
It follows that $f'$ has an injective inverse 
$f'^{-1} \in {\sf invfP}^{\sf (NP)}$, i.e., $f'$ is regular in 
${\sf cofP} = {\sf invfP}^{\sf (NP)}$.
This holds for all $f' \in {\sf cofP}$, hence 
${\sf cofP}$ is regular.
Thus by Theorem \ref{cofP_vs_P_NP}, ${\sf P} = {\sf NP}$.

\noindent $[\Rightarrow]$ \ If  ${\sf P} = {\sf NP}$ then
$\, {\sf cofP} = {\sf invfP}$, which in this case is also 
equal to ${\sf invfP} = {\sf invfP}^{\sf (NP)}$. 
 
\smallskip

If $e$ is an idempotent of ${\sf invfP}^{\sf (NP)}$ then $e$ 
is also a co-inverse of an element of {\sf fP}; indeed, 
$e \, \circ \, {\sf id}_{A^*} \, \circ \, e = e$. Hence 
$e \in {\sf cofP}$.
 \ \ \ $\Box$

\begin{pro} \label{RegIncofPvsinvfP}
 \ \ ${\sf P} \neq {\sf NP}$ \ iff \ ${\sf cofP}$ contains 
regular elements that do not belong to {\sf invfP}.
\end{pro}
{\bf Proof.} $[\Leftarrow]$ \ If ${\sf P} = {\sf NP}$ then 
${\sf cofP} = {\sf invfP}$ (by Theorem \ref{cofP_vs_P_NP}), 
so there are no (regular) elements in ${\sf cofP}$ that are not 
in {\sf invfP}.

\noindent $[\Rightarrow]$ \ If  ${\sf P} \neq {\sf NP}$ then we construct 
the following example of a regular element in ${\sf cofP}$ that 
does not belong to {\sf invfP}. Let $L \subset A^*$ be a {\sf coNP}-complete 
set.  Then ${\sf id}_L \in {\sf invfP}^{\sf (NP)}$ and it is regular (being 
an idempotent). Also, 
${\sf id}_L = {\sf id}_L \circ {\sf id}_{A^*} \circ {\sf id}_L$, hence 
${\sf id}_L$ is a co-inverse of an element of {\sf fP}. Thus,
${\sf id}_L \in {\sf cofP}$. But if ${\sf P} \neq {\sf NP}$ then
${\sf id}_L \not\in {\sf invfP}$, since $L$ is {\sf coNP}-complete.
 \ \ \ $\Box$

\bigskip

We saw that {\sf invfP} and ${\sf invfP}^{\sf (NP)}$ are finitely generated.
Hence, if {\sf cofP} were not finitely generated, this would imply that 
${\sf cofP} \neq {\sf invfP}$, which would imply that 
${\sf P} \neq {\sf NP}$ (by Theorem \ref{cofP_vs_P_NP}).
However, we will prove next that {\sf cofP} is finitely generated.  

We will first show that {\sf cofP} has a machine (or program) model, 
and that there is a corresponding evaluation function for all 
bounded-complexity functions in {\sf cofP}.  This is similar to the 
situation in {\sf fP}, {\sf invfP}, and ${\sf invfP}^{\sf (NP)}$. 

A {\em {\sf cofP}-program} is of the form $(v', w)$, where $v'$ is any
${\sf invfP}^{\sf (NP)}$-program (with built-in time-complexity and balance
function), and $w$ is any {\sf fP}-program (with built-in time-complexity 
and balance function). An ${\sf invfP}^{\sf (NP)}$-program describes a 
polynomially balanced polynomial-time injective Turing machine, with a fixed
{\sf NP}-oracle $N$ (where $N$ is {\sf NP}-complete).
The functions in {\sf fP} or ${\sf invfP}^{\sf (NP)}$ specified by programs
$w$ or $v'$ are denoted by $\phi_w$, respectively $\psi_{v'}$.
On input $y \in A^*$, the program $(v', w)$ is evaluated as follows: 

\smallskip

\noindent (1) \ \ \ $x = \psi_{v'}(y) \, $ is computed; 

\noindent (2) \ \ \ $\phi_w(x) \, $ is computed;

\noindent (3) \ \ \ if 
$\, \phi_w \circ \psi_{v'} \circ \phi_w(x) = \phi_w(x) \, $ and 
$ \, \psi_{v'} \circ \phi_w \circ  \psi_{v'}(y) = \psi_{v'}(y) \, $ then 
the output is $\psi_{v'}(y)$; 

 \ \ \ there is no output otherwise. 

\smallskip

\noindent The function specified by program $(v', w)$ on input $y$ is 
denoted by $\Phi_{(v', w)}$.
It is easy to see that $\Phi_{(v', w)}$ is a subinverse of $\phi_w$, and a
subfunction of $\psi_{v'}$. Thus, $\Phi_{(v', w)} \in {\sf cofP}$. 

Conversely, every function $h \in {\sf cofP}$ has an 
${\sf invfP}^{\sf (NP)}$-program, say $v'$;  and $h$ has an inverse 
$\phi_w \in {\sf fP}$ for some {\sf fP}-program $w$. Then $(v',w)$ is a
{\sf cofP}-program for $h$. 

\smallskip

Based on the {\sf cofP}-programs and a polynomial bound $q$, we can 
construct an evaluation function ${\sf evCo}_q^{(N)}$ such that  
$ \, {\sf evCo}_q^{(N)}({\sf code}(v') \, 11 \, {\sf code}(w) \, 11 \ y)$
$ \ = \ {\sf code}(v') \, 11 \, {\sf code}(w) \, 11 \ \Phi_{(v',w)}(y)$, 
 \ for $y$ and $(v', w)$ as above. This is similar to Section 3.2.

\begin{thm} \label{cofPFinGen}
 \ The monoid {\sf cofP} is finitely generated.
\end{thm}
{\bf Proof.} The proof is the same as for ${\sf invfP}^{\sf (NP)}$ (Prop.\ 
\ref{invfPNP_FinGen}), except that ${\sf evCo}_q$ (with programs $(v',w)$ as
above) is used instead of ${\sf injEv}_q^{(N)}$. 
 \ \ \ $\Box$

\bigskip

\noindent We saw that ${\sf P} \ne {\sf NP}$ iff {\sf cofP} is not regular. 
The set of elements of {\sf cofP} that are regular in
{\sf cofP} will be denoted by {\sf RegcofP}. We will see that it has
interesting properties.

\begin{pro} \label{RegMonCofP}
 \ (1) The set {\sf RegcofP} is a finitely generated inverse submonoid of
{\sf cofP}.  

\smallskip

\noindent (2) If ${\sf P} \ne {\sf NP}$ then 
$\ {\sf invfP} \, \subsetneqq \, {\sf RegcofP} \, \subsetneqq \, $
${\sf cofP} \, \subsetneqq \, {\sf invfP}^{\sf (NP)}$. 
 \ If ${\sf P} = {\sf NP}$ then 
$\, {\sf invfP} = {\sf invfP}^{\sf (NP)}$.    

\smallskip

\noindent (3) An element $f \in {\sf fP}$ has an inverse in {\sf RegcofP} 
iff $f$ is regular in {\sf fP}.
\end{pro}
{\bf Proof.} (1) For any $g \in {\sf RegcofP}$, let $g' \in {\sf RegcofP}$ 
be a mutual inverse of $g$. Multiplying $g g' g = g$ on the left and the 
right by $g^{-1}$ we obtain:
$\, {\sf id}_{{\sf Dom}(g)} \circ g' \circ {\sf id}_{{\sf Im}(g)} = g^{-1}$.
Hence $g^{-1} \in {\sf cofP}$, being a product of elements of {\sf cofP}
(since the idempotents ${\sf id}_{{\sf Dom}(g)}$,
${\sf id}_{{\sf Im}(g)}$ belong to {\sf cofP}).
Thus we proved:

\smallskip

 \ \ \  {\it For all $g \in {\sf RegcofP} : $
                $\, g^{-1} \in {\sf cofP}$.}

\smallskip

\noindent For any $g_1, g_2 \in {\sf RegcofP}$ we have therefore,
$g_1^{-1}, g_2^{-1} \in {\sf cofP}$, hence $g_2 g_1$ has
$g_1^{-1} g_2^{-1} \in {\sf cofP}$ as an inverse, so $g_2 g_1$ is regular 
in {\sf cofP}. This proves that {\sf RegcofP} is closed under composition.
Also, since all elements of {\sf RegcofP} are injective, {\sf RegcofP} is
an inverse monoid.

\smallskip

The proof of finite generation of {\sf RegcofP} is similar to the proof of
finite generation of {\sf cofP}.
We construct an evaluation function for the elements of {\sf RegcofP},
based on the following machine (or program) model for the elements of
{\sf RegcofP}. A {\sf RegcofP}-program is any pair $(u,v')$ of 
{\sf cofP}-programs, and the function $\Phi_{(u,v')}$ computed by this 
program is defined by

\smallskip

 \ \ \ $\Phi_{(u,v')}(x) \ = \ \psi_u(x)$ \ \ if
 \ \ $\psi_u \, \psi_{v'} \, \psi_u(x) = \psi_u(x)$;

\smallskip

\noindent and $\Phi_{(u,v')}(x)$ is undefined otherwise.
Since the relation $\psi_u \, \psi_{v'} \, \psi_u(x) = \psi_u(x)$
can be checked in ${\sf P}^{\sf (NP)}$, $\Phi_{(u,v')}$ belongs to 
${\sf invfP}^{\sf (NP)}$. And since $\psi_u \in {\sf cofP}$ has some 
inverse $f \in {\sf fP}$, and $\Phi_{(u,v')}$ is a subfunction of $\psi_u$,
$\Phi_{(u,v')}$ also has $f$ as an inverse; hence, 
$\Phi_{(u,v')} \in {\sf cofP}$. Moreover, since 
$\psi_u \, \psi_{v'} \, \psi_u(x) = \psi_u(x)$ for every 
$x \in {\sf Dom}(\Phi_{(u,v')})$, and $\Phi_{(u,v')} \subseteq \psi_u$,
we have $\Phi_{(u,v')} \, \psi_{v'} \, \Phi_{(u,v')} = \Phi_{(u,v')}$.
So $\Phi_{(u,v')}$ is regular in {\sf cofP}. 
Finally, every regular element of {\sf cofP} obviously has a program of the
form $(u,v')$ as above.

The rest of the proof of finite generation is very similar to the one for
{\sf cofP}.

\smallskip

\noindent (2) The second part of (2) was already proved in Theorem
\ref{cofP_vs_P_NP}. On the other hand, if ${\sf P} \ne {\sf NP}$ then
{\sf cofP} is not regular (by Theorem \ref{cofP_vs_P_NP}), so
${\sf RegcofP} \ne {\sf cofP} \ne {\sf invfP}^{\sf (NP)}$.
Also, {\sf cofP} contains regular elements that are not in {\sf invfP}, 
by Prop.\ \ref{RegIncofPvsinvfP}, hence ${\sf invfP} \ne {\sf RegcofP}$. 

\smallskip

\noindent (3) If $f$ is a regular element of {\sf fP} then $f$ has an
inverse in {\sf invfP}, by Coroll.\ \ref{injreg}.

In general, let $f$ be an element in {\sf fP} that has an inverse 
$f' \in {\sf RegcofP}$; we want to show that $f$ is regular in {\sf fP}.
By what we proved in (1) of the present proof, $f'^{-1} \in {\sf cofP}$.
We have $f'^{-1} g f'^{-1} = f'^{-1}$ for some $g \in {\sf fP}$, since 
$f'^{-1} \in {\sf cofP}$. Multiplying this on the left and the right by 
$f'$ yields  $\, {\sf id}_I \circ g \circ {\sf id}_D = f'$, where 
$D = {\sf Dom}(f')$ and $I = {\sf Im}(f')$. Since $f'$ is an inverse of $f$, 
Lemma \ref{regExt} implies that $g$ is an inverse of $f$. Since 
$g \in {\sf fP}$, it follows that $f$ is regular in {\sf fP}.
 \ \ \ $\Box$

\bigskip

We saw that the monoids {\sf fP}, {\sf invfP}, ${\sf invfP}^{\sf (NP)}$, 
and {\sf cofP} are finitely generated. 
For {\sf injfP} we ask similarly:  

\medskip

{\bf Question:} \ {\it  Is $\, {\sf injfP}$ finitely generated?}

\medskip

\noindent We do not know the answer.
If we could show that {\sf injfP} is {\sl not} finitely generated, then 
this would prove that ${\sf injfP} \neq {\sf invfP}$, i.e., injective 
one-way functions exist, 
hence ${\sf P} \neq {\sf UP}$, and hence ${\sf P} \neq {\sf NP}$. 

When we proved that {\sf fP}, {\sf invfP}, ${\sf invfP}^{\sf (NP)}$, and 
{\sf cofP} are finitely generated, we used machine (or program) models,
and evaluation functions.
It seems that {\sf injfP} does not have a program model, since injectiveness 
is a for-all property, that does not have finite witnesses in general. 
%%%%%%%%%%%%%%%%%%
To illustrate the difficulty of finding a program model for {\sf injfP}, 
here is an idea that does not work.  We saw that 
${\sf injfP} = {\sf fP} \cap {\sf invfP}^{\sf (NP)}$, and this suggests that
a function in {\sf injfP} can be specified by a pair $(u, w)$, where $w$ 
is an {\sf fP}-program (for $\phi_w \in {\sf fP}$), and $u$ is an 
${\sf invfP}^{\sf (NP)}$-program (for $\psi_u \in {\sf invfP}^{\sf (NP)}$); 
the program $(u, w)$ specifies the injective function $\psi_u \cap \phi_w$. 
But $\psi_u \cap \phi_w$ ranges over all of ${\sf invfP}^{\sf (NP)}$ and 
does not necessarily belong to {\sf injfP} (unless ${\sf injfP} =$ 
${\sf invfP}^{\sf (NP)}$, which would imply that ${\sf P} = {\sf NP}$). 
So, this approach towards proving finite generation does not work (unless 
one also proves that ${\sf P} = {\sf NP}$). 
It seems that {\sf injfP} is not finitely generated (but this will probably 
be very difficult to prove).

%%%%%%%%%%%%%%%%%%%%%%%%%%%%%%%%%%%%%%%%%%%%%%%%%%%
\section{Subgroups and group-inverses}

We first characterize the maximal subgroups of {\sf invfP}, 
{\sf cofP}, and ${\sf invfP}^{\sf (NP)}$.
Then we consider elements of {\sf fP} that have an inverse in such a
subgroup (i.e., a {\em group-inverse}).  Typically, elements of {\sf fP} do 
not have a group-inverse, but we will see that all elements of {\sf fP} are 
reduction-equivalent to elements that have group-inverses.

%%%%%%%%%%%
\subsection{Maximal subgroups}

We first characterize the idempotents of {\sf injfP}.

\begin{pro} \label{inj_Idempot_fP}
 \ For any $f \in {\sf injfP}$ we have:
$f$ is an {\em idempotent} \ iff \ $f = {\sf id}_Z$ for some set
$Z \subseteq A^*$ with $Z \in {\sf P}$.
Hence {\sf injfP} and {\sf invfP} have the same idempotents. 
\end{pro}
{\bf Proof.} If $f \in {\sf injfP}$ and $f = f \circ f$, then we compose on 
the left with $f^{-1}$ (which exists, since $f$ is injective, but $f^{-1}$ 
might not belong to {\sf injfP}); this yields 
${\sf id}_D = {\sf id}_D \circ f$, where $D = {\sf Dom}(f)$. Hence
for all $x \in D = {\sf Dom}(f)$ we have $x = f(x)$. Thus, 
$f = {\sf id}_D$.  Since $f \in {\sf fP}$, $D = {\sf Dom}(f) \in {\sf P}$. 

Conversely, if $Z \in {\sf P}$ then ${\sf id}_Z$ can be computed by a
deterministic polynomial-time Turing machine. Moreover,
${\sf id}_Z \in {\sf invfP}$ since ${\sf id}_Z$ is injective, and it is
regular (being an idempotent).
  \ \ \ $\Box$

\begin{pro} \label{groups_invfP_RM}
 \ The maximal subgroup of {\sf invfP} and {\sf injfP} are the same.

 \ Let ${\sf id}_Z$ be an idempotent of {\sf invfP}. 
The maximal subgroup of {\sf invfP} with identity ${\sf id}_Z$ consists of
the permutations of $Z$ that belong to {\sf invfP}.
In particular, the group of units of {\sf invfP} consists of all 
permutations of $A^*$ that belong to {\sf invfP}.  
\end{pro}
{\bf Proof.} If $f$ belongs to a subgroup of {\sf injfP} then $f$ is 
regular, hence $f \in {\sf invfP}$. 

Obviously, the permutations in {\sf invfP} with domain and image
$Z$ form a subgroup of {\sf invfP}.
Conversely, if $f \in {\sf invfP}$ belongs to the maximal subgroup with
identity ${\sf id}_Z$, then $f^{-1} f = f f^{-1} = {\sf id}_Z$; so the
domain and the image of $f$ are both $Z$. Since $f \in {\sf invfP}$, $f$ is
injective, and $f^{-1} \in {\sf invfP}$. Hence $f$ permutes $Z$.
 \ \ \ $\Box$

\begin{pro} \label{GroupsfPinvfP}
 \ Every maximal subgroup of {\sf invfP} 
is also a maximal subgroup of {\sf fP}.
The group of units of {\sf fP} is the same as the group of units of 
{\sf invfP}.
\end{pro}
{\bf Proof.} A maximal subgroup of {\sf fP} is a regular $\cal H$-class.
For {\sf fP} 
the characterization of the
$\cal L$- and the $\cal R$-relation (see Prop.\ 2.1 in \cite{s1f})
implies that every element $f$ of a maximal subgroup of {\sf fP}
with unit ${\sf id}_Z$ satisfies ${\sf Im}(f) = Z = {\sf Dom}(f)$ and $f$ 
is injective (since $f$ has the same partition as ${\sf id}_Z$). Thus $f$ 
belongs to {\sf invfP}.
 \ \ \ $\Box$

\medskip

\noindent The converse of Prop.\ \ref{GroupsfPinvfP} is of course not true;
e.g., {\sf fP} contains non-injective idempotents, hence subgroups that are
not contained in {\sf invfP}.

\bigskip

Let us now look at the idempotents and subgroups of ${\sf invfP}^{\sf (NP)}$
and of {\sf cofP}. We saw that ${\sf invfP}^{\sf (NP)}$ and 
{\sf cofP} have the same idempotents (Prop.\ \ref{cofPvsfPNP}).  
And we saw that ${\sf invfP}^{\sf (NP)}$ is an inverse monoid (Cor.\ 
\ref{invfPoracle_inv}); on the other hand, {\sf cofP} is not regular,
unless ${\sf P} = {\sf NP}$ (by Theorem \ref{cofP_vs_P_NP}).

\begin{lem} \label{L_vs_idL}
 \ For every $L \subset A^*$ we have: \ $L \in {\sf P}^{\sf (NP)}$ \ iff 
 \ ${\sf id}_L \in {\sf invfP}^{\sf (NP)}$. 
\end{lem}
{\bf Proof.} $[\Rightarrow]$ If $L \in {\sf P}^{\sf (NP)}$ then ${\sf id}_L$ 
can be computed by deterministic polynomial-time Turing machine with 
{\sf NP}-oracle. 
Moreover, the inverse of ${\sf id}_L$ (which is just ${\sf id}_L$ itself) 
also has that property, hence by Prop.\ \ref{injVSdetandrev}, ${\sf id}_L$ 
can be computed by an injective polynomial-time Turing machine with 
{\sf NP}-oracle. Hence, by the definition of ${\sf invfP}^{\sf (NP)}$ we have
${\sf id}_L \in {\sf invfP}^{\sf (NP)}$.

$[\Leftarrow]$ From a polynomial-time injective Turing machine with
{\sf NP}-oracle, computing ${\sf id}_L$, one immediately obtains a
deterministic polynomial-time Turing machine with {\sf NP}-oracle, accepting 
$L$.  
 \ \ \ $\Box$

\begin{pro} \label{idemp_cofP}
 \ Every idempotent of $\, {\sf invfP}^{\sf (NP)}$ (and {\sf cofP}) is 
of the form ${\sf id}_L$, where $L \in {\sf P}^{\sf (NP)}$, 
 \ $L \subseteq A^*$.
Moreover, every $L \in {\sf P}^{\sf (NP)}$ is accepted by an {\em injective} 
polynomial-time Turing machine with {\sf NP}-oracle.  
\end{pro}
{\bf Proof.} Let $f = f \circ f$ be an idempotent of 
${\sf invfP}^{\sf (NP)}$.
Multiplying in the left by $f^{-1}$ yields ${\sf id}_D = {\sf id}_D \circ f$,
where $D = {\sf Dom}(f)$. Hence, $f = {\sf id}_D$. Since $f \in$ 
${\sf invfP}^{\sf (NP)}$ we have $D \in {\sf P}^{\sf (NP)}$ (by Lemma 
\ref{L_vs_idL}).

Conversely, if $Z \in {\sf P}^{\sf (NP)}$ then ${\sf id}_Z \in$
${\sf invfP}^{\sf (NP)}$ (by Lemma \ref{L_vs_idL}).
From this we obtain an injective polynomial-time Turing machine with 
{\sf NP}-oracle, accepting $Z$.
 \ \ \ $\Box$

\begin{pro} \label{groupsCofp} 
 \  (1) The maximal subgroup of ${\sf invfP}^{\sf (NP)}$ with idempotent 
${\sf id}_L$ (where $L \in {\sf P}^{\sf (NP)}$) consists of all the 
permutations of $L$ that belong to ${\sf invfP}^{\sf (NP)}$.
In particular, the group of units of ${\sf invfP}^{\sf (NP)}$ consists of 
all permutations of $A^*$ that belong to ${\sf invfP}^{\sf (NP)}$.         

\smallskip

\noindent 
(2) The maximal subgroup of {\sf cofP} with idempotent ${\sf id}_L$ (where
$L \in {\sf P}^{\sf (NP)}$) consists of all permutations $g$ of $L$ such 
that $g = f \, \cap \, (L \times L)$ and 
$g^{-1} = h \, \cap \, (L \times L)$ for some $f, h \in {\sf fP}$. 

Hence, if $L \in {\sf P}$ then the maximal subgroup of {\sf cofP} with
idempotent ${\sf id}_L$ is a subgroup of {\sf invfP}. In particular, the 
group of units of {\sf cofP} is the same as the group of units of
{\sf invfP}.
\end{pro}
{\bf Proof.} (1) By injectiveness of the elements of 
${\sf invfP}^{\sf (NP)}$, every element of the maximal subgroup of 
${\sf invfP}^{\sf (NP)}$ with idempotent ${\sf id}_L$ is a permutation of 
$L$. Moreover, for every $g \in {\sf invfP}^{\sf (NP)}$ we have 
$g^{-1} \in {\sf invfP}^{\sf (NP)}$ (by Prop.\ \ref{Prop_invfPNP}); hence 
the permutations of $L$ in ${\sf invfP}^{\sf (NP)}$ form a subgroup of
${\sf invfP}^{\sf (NP)}$.  

\smallskip

(2) If $g \in {\sf invfP}^{\sf (NP)}$ belongs to a maximal subgroup of 
{\sf cofP} with idempotent ${\sf id}_L$, then both $g$ and $g^{-1}$ are
permutations of $L$. If $g$ has an inverse $h \in {\sf fP}$ then multiplying 
$g h g = g$ on the left and on the right by $g^{-1}$ yields 
${\sf id}_L \circ h \circ {\sf id}_L = g^{-1}$; 
moreover, ${\sf id}_L \circ h \circ {\sf id}_L = h \, \cap \, (L \times L)$.
Similarly, for $g^{-1} \in {\sf coP}$ we obtain 
$f \, \cap \, (L \times L) = g$ for some $f \in {\sf fP}$. 

Conversely, let $g$ be a permutation of $L \in {\sf P}^{\sf (NP)}$ such that 
$g = f \, \cap \, (L \times L)$ and $g^{-1} = h \, \cap \, (L \times L)$ for 
some $f, h \in {\sf fP}$. Then $g, g^{-1}$ belong to 
${\sf invfP}^{\sf (NP)}$ since $L \in {\sf P}^{\sf (NP)}$ and 
$f, h \in {\sf fP}$. 
Moreover, $g$ and $g^{-1}$ belong to {\sf cofP}; indeed, the
inverses $f \, \cap \, (L \times L)$ and $h \, \cap \, (L \times L)$ can be
extended to inverses $h$, respectively $f$, in {\sf fP} (by Lemma 
\ref{regExt}).
Since $g$ is a permutation of $L$, it follows now that $g$ belongs to a 
subgroup of {\sf cofP}.
 \ \ \ $\Box$

\bigskip

We saw in Prop.\ \ref{RegIncofPvsinvfP} that unless ${\sf P} = {\sf NP}$, 
{\sf cofP} contains regular elements that are not in {\sf invfP}.
The next Proposition describes the regular elements of {\sf cofP}.

\begin{pro} \label{regINcofP}
 \ Every element of {\sf cofP} that is regular in {\sf cofP} has 
the form $f \, \cap \, (K \times H)$, for some $f \in {\sf fP}$ and 
$K,H \in {\sf P}^{\sf (NP)}$. 

Conversely, suppose $f, h \in {\sf fP}$ and $K, H \in {\sf P}^{\sf (NP)}$
are such that $f \, \cap \, (K \times H)$ and $h \, \cap \, (H \times K)$
are injective and mutual inverses.
Then $f \, \cap \, (K \times H)$, $h \, \cap \, (H \times K)$ $\in$
${\sf cofP}$, and they are regular in {\sf cofP}.  
\end{pro} 
{\bf Proof.} If $g$ is regular in {\sf cofP}
then $g^{-1} \in {\sf cofP}$; hence, $g^{-1} \, f \, g^{-1} = g^{-1}$
for some $f \in {\sf fP}$. Hence by multiplying on the left and the right by
$g$ we obtain: $f \, \cap \, (K \times H) = g$, where $K = {\sf Im}(g)$
and $H = {\sf Dom}(g)$.
 
For the converse, $f \, \cap \, (K \times H) \ = \ $
$(f \, \cap \, (K \times H)) \circ (h \, \cap \, (H \times K))$ $\circ$
$(f \, \cap \, (K \times H))$ $ \ = \ $
$(f \, \cap \, (K \times H)) \circ h \circ (f \, \cap \, (K \times H))$;
the latter holds by Lemma \ref{regExt}.
Hence, $f \, \cap \, (K \times H)$ has an inverse in {\sf fP}, so 
$f \, \cap \, (K \times H)$ belongs to {\sf cofP}.
Similarly, $h \, \cap \, (H \times K)$ belongs to {\sf cofP}.
Since they are mutual inverses, they are regular in {\sf cofP}.
 \ \ \ $\Box$

%%%%%%%%%%%
\subsection{Group inverses of elements of {\sf fP} }

After noticing that every regular element of {\sf fP}
has an injective inverse, belonging to an inverse submonoid of {\sf fP}, 
we wonder whether we can go even further: Does every regular element in 
{\sf fP} have an inverse in a {\em subgroup} of {\sf fP}?
This is of course not the case; for immediate counter-examples consider the 
functions that are total but not surjective, or surjective but not total. 
Even within the subsemigroup of non-total non-surjective functions of 
{\sf fP} there are counter-examples (due to the polynomial balance and time
requirements).
Nevertheless, we will find that every regular element of {\sf fP} is
equivalent, with respect to inversive reduction, to a regular element of 
{\sf fP} that has an inverse in a subgroup of {\sf fP}.

We will also investigate elements of {\sf fP} that are possibly non-regular,
but that have an inverse in a subgroup of ${\sf invfP}^{\sf (NP)}$. Again, we 
will find that every element of {\sf fP} is equivalent, with respect to
inversive reduction, to an element of {\sf fP} that has an inverse in a 
subgroup of ${\sf invfP}^{\sf (NP)}$.
In particular, there are elements of {\sf fP} that are complete with respect 
to inversive reduction and that have an inverse in a subgroup of 
${\sf invfP}^{\sf (NP)}$.

\bigskip

\noindent Some definitions about inversive reductions (from \cite{s1f}) 
between functions $f_1, f_2$:

\smallskip

\noindent $\bullet$ \ \ $f_1$ is {\em simulated} by $f_2$ (denoted by 
$f_1 \preccurlyeq f_2$) \ iff \ there exist $\beta, \alpha \in {\sf fP}$ 
such that $f_1 = \beta \circ f_2 \circ \alpha$.

\smallskip

\noindent $\bullet$ \ \ $f_1$ {\em reduces inversively} to $f_2$ 
(notation, $f_1 \leqslant_{\sf inv} f_2$) \ iff 

 \ \ \ \ \ {\small (1)} \ $f_1 \preccurlyeq f_2$ \ and 

 \ \ \ \ \ {\small (2)} \ for every inverse $f_2'$ of $f_2$ there exists 
    an inverse $f_1'$ of $f_1$ such that $f_1' \preccurlyeq f_2'$ .

\smallskip

\noindent $\bullet$ \ \ $f, g \in {\sf fP}$ are {\em equivalent via 
inversive reduction} 
 \ iff \ $f \leqslant_{\sf inv} g \,$ and $\, g \leqslant_{\sf inv} f$.

%%%%%%%%%%%%%%%%%% Theorem 
\begin{thm} \label{grpinvAllfP}
 \ Every function $f \in {\sf fP}$ is equivalent, via inversive reduction, 
to a function $f_0 \in {\sf fP}$ such that $f_0$ has an inverse in a 
{\em subgroup} of ${\sf invfP}^{\sf (NP)}$. 

In particular, $f_0$ has an inverse in the group of permutations of 
$\  0 \, {\sf Dom}(f) \, \cup \, 1 \, {\sf Im}(f)$, in 
${\sf invfP}^{\sf (NP)}$. 
Moreover, $f_0$ has an inverse in the group of units of 
${\sf invfP}^{\sf (NP)}$.

If $f$ (and hence $f_0$) is {\em regular}, then $f_0$ has an inverse in 
the group of permutations of 
$\, 0 \, {\sf Dom}(f) \, \cup \, 1 \, {\sf Im}(f)$, in {\sf invfP}.
Moreover, $f_0$ has an inverse in the group of units of {\sf invfP}. 
\end{thm}
{\bf Remarks.} Since $f$ and $f_0$ are equivalent via reduction it follows 
immediately that $\, f_0 \in {\sf fP}$ iff $f \in {\sf fP}$; and $f$ is 
regular iff $f_0$ is regular. 
And by Theorem \ref{cofP_vs_P_NP}, if $f$ is regular then $f_0$ has an 
inverse in the inverse monoid {\sf invfP}.  

\medskip

\noindent {\bf Proof of Theorem \ref{grpinvAllfP}.}
(1) With $f$ we associate $f_0$ defined by

\smallskip

 \ \ \ ${\sf Dom}(f_0) \ = \ 0 \ {\sf Dom}(f)$, \ and 

\smallskip

 \ \ \ $f_0(0 \, x) = 1 \, f(x)$, \ for all $x \in {\sf Dom}(f)$. 

\smallskip

\noindent Similarly, we define $f_1$ by $\, f_1(1 \, x) = 0 \, f(x)$.

To show that $f$ and $f_0$ simulate each other (and similarly for $f_1$), 
we introduce the functions $\pi_a$ and $\pi_a'$ for each $a \in A$; they 
are defined for all $z \in A^*$ by

\smallskip

 \ \ \ $\pi_a(z) = a z$, \ \ and \ \ $\pi_a'(az) =z$, with 
$\, {\sf Dom}(\pi_a') = aA^*$.

\smallskip

\noindent Then we have:

\smallskip

 \ \ \ $f_0 = \pi_1 \circ f \circ \pi_0'$, \ \ and 
 \ \ $f = \pi_1' \circ f_0 \circ \pi_0$. 

\smallskip

\noindent Hence, $f$ and $f_0$ simulate each other. 

Let us show that we have an {\em inversive} reduction of $f_0$ to $f$.  
If $f$ has an inverse $f'$, let us define $f'_1$ by
$f'_1(1 \, y) = 0 \, f'(y)$ for all $y \in {\sf Dom}(f')$. Then $f'_1$ is 
an inverse of $f_0$. And $f'_1 = \pi_0 \circ f' \circ \pi_1'$, so 
$f_1'$ is simulated by $f'$. 

Conversely, let us show that there is an inversive reduction of $f$ to 
$f_0$.  Let $g$ be any inverse of $f_0$, i.e., 
$f_0 g f_0(x) = f_0(0 \, x)$,  all $x \in {\sf Dom}(f)$.  
Then $h = g \, \cap \, 1 A^* \times 0 A^*$ can be simulated by $g$ (indeed,
$h(.) = {\sf id}_{0 A^*} \circ g \circ {\sf id}_{1  A^*}(.)$), and 
$h$ also satisfies $f_0 h f_0 = f_0$. Moreover, 
$h \subset 1 A^* \times 0 A^*$ implies that $h = k_1$ for some function
$k$. Then we have $f k f = f$; indeed, for all $x \in {\sf Dom}(f)$ we have
$f k f(x) = \pi_1' \, f_0 \, k_1 \, f_0(0x) = \pi_1' \, f_0(0x) = f(x)$.  
We saw that $k$ is simulated by $k_1$ (indeed, 
$k = \pi_0' \circ k_1 \circ \pi_1$). 
Thus there exists an inverse (namely $k$) of $f$ that is simulated by $g$.

This completes the proof that $f$ is equivalent to $f_0$ via inversive 
reduction. It follows immediately that $f_0 \in {\sf fP} \, $ iff 
$\, f \in {\sf fP}$, and that $f$ is regular iff $f_0$ is regular. 

\smallskip

\noindent (2) \ Let $f'$ be any mutual inverse of $f$ with 
$f' \in {\sf cofP}$, or with $f' \in {\sf invfP}$ if $f$ is regular.
Let us now extend $f_1'$ to a group element.
Recall that by definition, $f'_1(1 \, y) = 0 \, f'(y)$ for all 
$y \in {\sf Dom}(f')$, 
$\, {\sf Dom}(f_1') = 1 \, {\sf Dom}(f') = 1 \, {\sf Im}(f)$, and 
${\sf Im}(f_1') = 0 \, {\sf Im}(f)$. 
So, ${\sf Dom}(f_1') \, \cap \,{\sf Im}(f_1') = \varnothing$, and 
$f_1'$ only has orbits of length 1.  We extend $f_1'$ to a 
permutation $F'$ of $\, {\sf Dom}(f_1') \cup {\sf Im}(f_1')$, defined by

\smallskip

 \ \ \ $F'|_{{\sf Dom}(f_1')} = f_1'$  \ \ \ and
 \ \ \ $F'|_{{\sf Im}(f_1')} = f_1'^{-1}$.

\smallskip

\noindent Hence, $F' = f_1' \cup f_1'^{-1}$.
It follows that $F' \in {\sf invfP}^{\sf (NP)}$, since both $f_1'$ and 
$f_1'^{-1}$ belong to ${\sf invfP}^{\sf (NP)}$. 

Also, $F'$ is injective and $\, F' \circ F' = {\sf id}_Z$, where 
$Z = {\sf Dom}(f_1') \cup {\sf Dom}(f_1'^{-1})$;  so $F' = F'^{-1}$ belongs 
to a two-element group.  Clearly, ${\sf Dom}(F') = {\sf Im}(F') = Z$.

Moreover, $F'$ is an inverse of $f_0$. Indeed, for every 
$x \in {\sf Dom}(f)$ we have $\, f_0 \, F' \, f_0(0 x) \, = $ 
$f_0 \, F'(1 f(x)) \, = \, f_0 \, f_1'(1 f(x))$, since 
$F' = f_1' \cup f_1'^{-1}$ and $f_1'^{-1}$ is undefined on $1 A^*$;
hence, $f_0 \, F' \, f_0(0 x) \, = \, f_0(0 x)$.

Finally, let us show that if $f' \in {\sf invfP}$ then $F' \in {\sf invfP}$.
Indeed, when $f' \in {\sf invfP}$ then $f'^{-1} \in {\sf invfP}$, from which
it follows (by disjointness of the domains and disjointness of the images) 
that $f_1' \cup f_1'^{-1} \in$ {\sf invfP}.

\smallskip

\noindent (3) \ We can extend $F'$ to a permutation of all of $A^*$ 
(belonging to the group of units of ${\sf invfP}^{\sf (NP)}$): We simply 
define $F'$ as a permutation of ${\sf Dom}(f_1') \cup {\sf Im}(f_1')$, as 
above, and then extend $F'$ to the identity function on
$A^* - ({\sf Dom}(f_1') \cup {\sf Im}(f_1'))$.  Since 
$f_1', \, {\sf id}_{{\sf Dom}(f_1')}$, and ${\sf id}_{{\sf Im}(f_1')}$
belong to ${\sf invfP}^{\sf (NP)}$, it follows that this extension of $F'$ 
to all of $A^*$ belongs to ${\sf invfP}^{\sf (NP)}$. 

Finally, if $f$ is regular then this extension of $F'$ belongs to 
{\sf invfP}; indeed, in that case, ${\sf Dom}(f_1')$ and ${\sf Im}(f_1')$ 
belong to {\sf P}. 
 \ \ \ $\Box$

\begin{cor} \label{ExistsGroupExtComplete}
 \ There exists a function that is {\em complete} in {\sf fP} for inversive
reduction, and that has an inverse in the group of units of  
${\sf invfP}^{\sf (NP)}$.
\end{cor}
{\bf Proof.} In \cite{s1f} we saw examples of functions $F \in {\sf fP}$
that are complete with respect to inversive reduction. Then $F_0$, defined
by $F_0(0 x) = 1 \, F(x)$, is also complete (since $F$ and $F_0$ are
equivalent via inversive reduction); and $F_0$ has an inverse in the
group of units of ${\sf invfP}^{\sf (NP)}$ (as we saw in the proof of 
Theorem \ref{grpinvAllfP}).
 \ \ \ $\Box$

\bigskip

\noindent {\bf Remark.} The group-inverses constructed above belong to  
${\sf invfP}^{\sf (NP)}$, but not necessarily to {\sf cofP}.
If there were a complete function in {\sf fP} that has an inverse in the
group of units of {\sf cofP} then ${\sf P} = {\sf NP}$ (since the 
group of units of {\sf cofP} is the same as the group of units of
{\sf invfP}).

%%%%%%%%%%%%%%%%%%%%%%%%%%%%%%%%%%%%%%%%%%%%%%%%%%%%%%

\section{Appendix: Simple facts about inverses (related to Section 2)}

We present some additional properties of {\sf invfP} that are not used in
the rest of the paper. 

For $f \in {\sf fP}$ we define the {\em right fixator} of $f$ by
$\, {\sf RFix}(f) = \{\alpha \in {\sf fP} :  f \circ \alpha = f \}$; in
other words, $\alpha \in {\sf RFix}(f)$ iff $\alpha$ is a right-identity of
$f$.  Similarly, the {\em left fixator} of $f$ is
 ${\sf LFix}(f)  = \{\alpha \in {\sf fP} : \alpha \circ f = f\}$.

We observe that ${\sf RFix}(f) \cap {\sf invfP}$ is an inverse monoid;
i.e., $\alpha \in {\sf RFix}(f) \cap {\sf invfP}$ implies
$\alpha^{-1} \in {\sf RFix}(f) \cap {\sf invfP}$.
Indeed, if $\alpha \in {\sf invfP}$ satisfies $f \alpha = f$ then
$f \alpha \alpha^{-1} = f \alpha^{-1}$; moreover, $f \alpha = f$
implies that ${\sf Dom}(f) \subseteq {\sf Im}(\alpha)$, hence
$\alpha \alpha^{-1} = {\sf id}_{{\sf Im}(\alpha)}$ acts as a
right-identity on $f$; hence $f = f \alpha \alpha^{-1}$.
Thus, $f = f \alpha^{-1}$.

For every $\alpha \in {\sf RFix}(f)$ we have:
 \ $\alpha({\sf Dom}(f)) \subseteq {\sf Dom}(f)$,
 \ $\alpha|_{{\sf Dom}(f)} \in {\sf RFix}(f)$, and
$\, \alpha(A^* - {\sf Dom}(f)) \subseteq  A^* - {\sf Dom}(f)$.

Every representative choice function of $f$ belongs to ${\sf RFix}(f)$, since
any right fixator $\alpha \in {\sf RFix}(f)$ maps every ${\sf mod}_f$-class
into itself.  But the converse is not true, since a right fixator $\alpha$
does not necessarily map a ${\sf mod}_f$-class to a single element of the
${\sf mod}_f$-class (as a representative choice function does).

\begin{pro} \label{prop_inverses}
 \ Let $f_1', f_2' \in {\sf invfP}$ be inverses of the same element
$f \in {\sf fP}$, such that ${\sf Dom}(f_i') = {\sf Im}(f)$ ($i = 1,2$).
Then:

\smallskip

\noindent {\bf (1)} \ \ \ $f_1'^{-1} \equiv_{\cal R} f_2'^{-1}$
$\equiv_{\cal R} {\sf id}_{{\sf Im}(f)} \, $ (for the ${\cal R}$-relation
of {\sf invfP}), \ and

\smallskip

\hspace{0.11in} $\, {\sf id}_{{\sf Im}(f)} \equiv_{\cal R} f \,$ (for
the ${\cal R}$-relation of {\sf fP}); and

\smallskip

\hspace{0.11in}  $\, f_1' \equiv_{\cal L} f_2' \, $ (for the
${\cal L}$-relation of {\sf invfP}).

\smallskip

\hspace{0.11in}  If $\, f_1' \equiv_{\cal L} f_2' \, $ and
$\, f_1' \equiv_{\cal R} f_2' \, $ (for the ${\cal L}$- and
${\cal R}$-relations of {\sf fP}), then $\, f_1' = f_2'$.

\smallskip

\noindent {\bf (2)} \ \ \ $f_i' f \in {\sf RFix}(f)$,

\smallskip

\hspace{0.11in} $f_2' f_1'^{-1}$ is a bijection from the choice set
${\sf Im}(f_1')$ onto the choice set ${\sf Im}(f_2')$;

\smallskip

\hspace{0.11in}  $\, f_i' f_i'^{-1} = {\sf id}_{{\sf Im}(f_i')}$.

\smallskip

\hspace{0.11in}  $\, f_i'^{-1} f_i' = {\sf id}_{{\sf Im}(f)}$.

\smallskip

\noindent {\bf (3)} \ \ \ $f_2' \, f \, f_1' = f_2'$.

\smallskip

\noindent {\bf (4)} \ \ \ If $\alpha \in {\sf RFix}(f)$, and if
$f' \in {\sf fP}$ is an inverse of $f \in {\sf fP}$ such that
${\sf Dom}(f') = {\sf Im}(f)$, then
$\alpha \circ f'$ is also an inverse of $f$ such that
${\sf Dom}(\alpha f') = {\sf Im}(f)$.
So, the monoid ${\sf RFix}(f)$ acts on the left on the set

\smallskip

\hspace{0.4in}  $\{ g' \in {\sf invfP} : \, g'$ {\rm is an inverse of $f$
and} ${\sf Dom}(g') = {\sf Im}(f)\}$.

\smallskip

\noindent Moreover, this action is transitive, and the action of
$\, \{\alpha|_{{\sf Dom}(f)} : \alpha \in {\sf RFix}(f)\} \,$ is
faithful.
\end{pro}
{\bf Proof.} (1) We saw that $f_i'^{-1} = f|_{{\sf Im}(f_i')}$, hence
$f_i'^{-1} = f \circ {\sf id}_{{\sf Im}(f_i')}$, where
${\sf Im}(f_i') \in {\sf P}$ (since $f_i'$ is regular); hence
${\sf id}_{{\sf Im}(f_i')} \in {\sf fP}$ and $f_i'^{-1} \le_{\cal R} f$.

Also, $f_i' = f_i' f f_i' \equiv_{\cal L} f f_i' = {\sf id}_{{\sf Im}(f)}$.
We also have $f_i'^{-1} \equiv_{\cal R} {\sf id}_{{\sf Im}(f)}$; indeed,
$f_i'^{-1}$ maps ${\sf Im}(f_i')$ bijectively onto ${\sf Im}(f)$, so
$f_i'^{-1} = {\sf id}_{{\sf Im}(f)} \circ f_i'^{-1} \le_{\cal R}$
${\sf id}_{{\sf Im}(f)}$; and ${\sf id}_{{\sf Im}(f)} = f_i'^{-1} f_i'$
$\le_{\cal R} f_i'^{-1}$.
Moreover, $f \equiv_{\cal R} {\sf id}_{{\sf Im}(f)}$, hence
$f_i'^{-1} \equiv_{\cal R} f$.

The $\, \equiv_{\cal L}$-equivalence then follows, since in any inverse
monoid, $x' \equiv_{\cal R} y'$ implies $x \equiv_{\cal L} y$ (where $z'$
is the inverse of $z$ for any element $z$).

The fact that $f_1' \equiv_{\cal L} f_2'$ and $f_1' \equiv_{\cal R} f_2'$
imply $f_1' = f_2'$ is well-known in semigroup theory (see e.g.\ p.\ 26 in 
\cite{CliffPres}).

(2) Obviously, $f \, f_i' \, f = f$.
The rest is straightforward.

(3) We know that $f_1'$ and $f_2'$ have the same domain, namely
${\sf Im}(f)$. And $f \, f_1'|_{{\sf Im}(f)}  = {\sf id}_{{\sf Im}(f)}$,
hence $f_2' \, f \, f_1' = f_2' \, {\sf id}_{{\sf Im}(f)} = f_2'$ (the
latter again holds since ${\sf Dom}(f_i') = {\sf Im}(f)$).

(4) If $f \alpha = f$ then $f \alpha f' f = f f' f = f$,
so $\alpha f'$ is also an inverse of $f$.
And ${\sf Dom}(\alpha f') \subseteq {\sf Dom}(f') = {\sf Im}(f)$;
since $\alpha f'$ is an inverse of $f$ we also have
${\sf Im}(f) \subseteq {\sf Dom}(\alpha f')$; hence
${\sf Dom}(\alpha f') = {\sf Im}(f)$.
It follows that $\alpha f' \in {\sf invfP}$, so ${\sf RFix}(f)$ acts on
the given set on the left.

Transitivity follows from $(f_2' f) \, f_1' = f_2'$ (proved in (3)),
where $f_2' f \in {\sf RFix}(f)$ by (2).

Proof of faithfulness: Note that ${\sf Dom}(f) \subseteq {\sf Dom}(\alpha)$
for all $\alpha \in {\sf RFix}(f)$. 
If $\alpha_1|_{{\sf Dom}(f)} \neq \alpha_2|_{{\sf Dom}(f)}$ then there 
exists $x_0 \in {\sf Dom}(f)$ such that 
$\alpha_1(x_0) \neq \alpha_2(x_0)$.
There exists an inverse $f_0' \in {\sf invfP}$ of $f$ with
${\sf Dom}(f_0') = {\sf Im}(f)$ such that $f_0'(f(x_0)) = x_0$; indeed, we
can start with any inverse $f_0'$ with ${\sf Dom}(f_0') = {\sf Im}(f)$,
and if $f_0'(f(x_0)) \neq x_0$ we can redefine $f_0'$ on $f(x_0)$; this
does not change ${\sf Dom}(f_0')$ (it just changes the choice set of 
$f_0'$). Now, 
$\alpha_1 f_0' f(x_0) = \alpha_1(x_0) \neq  \alpha_2(x_0)$ $=$
$\alpha_2 f_0' f(x_0)$, so $\alpha_1 f_0' \neq \alpha_2 f_0'$.
 \ \ \ $\Box$

\begin{pro}
 \ If $f \in {\sf fP}$ is regular, and if
$\, {\sf Dom}(f \circ \alpha) = {\sf Dom}(f)$,  then we have:

\smallskip

 \ \ \ $\alpha \in {\sf RFix}(f)$ \ \ iff \ \ for every mutual inverse
$f'$ of $f$:  $\alpha \circ f'$ is also a mutual inverse of $f$.

\smallskip

\noindent (If $f$ is not regular, the Proposition holds for inverses in
${\sf invfP}^{\sf (NP)}$.)
\end{pro}
{\bf Proof.} $[\Rightarrow]$ \ If $\alpha \in {\sf RFix}(f)$ and $f'$ is a
mutual inverse of $f$ then $f \, \alpha f' \, f = f f' f = f$;
the first equality holds because $\alpha \in {\sf RFix}(f)$. And 
$\alpha f'\, f \, \alpha f' = \alpha f'\, f \, f' = \alpha f'$ (again, 
using $f \alpha = f$). So $\alpha f'$ is a mutual inverse of $f$.

\noindent $[\Leftarrow]$ \ If $\alpha \not\in {\sf RFix}(f)$ then
$f \alpha(x_0) \neq f(x_0)$ for some $x_0 \in A^*$.
If $x_0 \in {\sf Dom}(f)$, let $f' \in {\sf fP}$ be a mutual inverse of
$f$ such that $f' f(x_0) = x_0$; such a mutual inverse exists (if
$f'f(x_0) \neq x_0$ we can redefine $f'$ on the one element $f(x_0)$).
Then $f \, \alpha f' \, f(x_0) = f \alpha(x_0) \neq f(x_0)$, so
$f \, \alpha f' \, f \neq f$. So there exists a mutual inverse $f'$ such
that $\alpha f'$ is not an inverse of $f$.
 \ \ \ $\Box$

\medskip

\noindent From the above we conclude that finding an inverse $f'$ for $f$
with ${\sf Dom}(f') = {\sf  Im}(f)$ can be broken up into two steps:
(1) find a maximal injective subfunction $f'^{-1}$ of $f$ (being a maximal
injective subfunction is equivalent to
${\sf Im}(f'^{-1}) = {\sf  Im}(f)$); (2) find the inverse $f'$ of
$f'^{-1}$. (If $f$ is injective, step (1) is skipped.)

\begin{pro}
 \ When $f \in {\sf fP}$ is regular and if
$\, {\sf Im}(\beta \circ f) = {\sf Im}(f)$ we have:

\smallskip

 \ \ \ $\beta \in {\sf LFix}(f)$ \ \ iff \ \ for every mutual inverse $f'$
of $f$: \ $f' \beta$ is also a mutual inverse of $f$.

\smallskip

\noindent (If $f$ is not regular, the Proposition holds for inverses in
${\sf invfP}^{\sf (NP)}$.)
\end{pro}
{\bf Proof.} $[\Rightarrow]$ \ If $\beta \in {\sf LFix}(f)$ and $f'$ is a
mutual inverse of $f$ then $f \, f' \beta \, f =  f f' f = f$; and
$f' \beta \, f \, f' \beta = f' \, f \, f' \beta = f' \beta$ (using 
$\beta f = f$). 

\noindent $[\Leftarrow]$ If $\beta \not\in {\sf LFix}(f)$ then
$\beta f(x_0) \neq f(x_0)$ for some $x_0 \in A^*$.
If $\beta f(x_0) \in {\sf Im}(f)$ then
$f \, f' \beta \, f(x_0) = \beta \, f(x_0)$ since
$f f'|_{{\sf Im}(f)} = {\sf id}_{{\sf Im}(f)}$; and
$\beta \, f(x_0) \neq f(x_0)$, hence $f \, f' \beta \, f \neq f$.
  \ \ \ $\Box$

\bigskip

\bigskip

%%%%%%%%%%%%%%%%%%%%%%%%%%%%%%%%%%%%%%%%%%%%%%%%%%%%%%%%%%%%%%%%
%% \small
%%%%%%%%%%%%%%%%%%%%%%%%%%%%%%%%%%%%%%%%%%%%%%%%%%%%

\end{document}